\documentclass{article}

\usepackage[final,nonatbib]{neurips_2020}

\usepackage[sort&compress]{natbib}
\setcitestyle{numbers,square,comma}

\usepackage[utf8]{inputenc} 
\usepackage[T1]{fontenc}    
\usepackage{hyperref}       
\usepackage{url}            
\usepackage{booktabs}       
\usepackage{amsfonts}       
\usepackage{nicefrac}       
\usepackage{microtype}      
\usepackage{hyperref}
\hypersetup{breaklinks=false,
bookmarks=false,
colorlinks=true,
citecolor=blue,
urlcolor=magenta,
linkcolor=magenta
}

\usepackage{tabu}

\usepackage{macros}
\usepackage{xfrac}

\usepackage[font=footnotesize]{caption}
\pgfplotsset{compat=1.14}
\usetikzlibrary{arrows,shapes,trees}
\usetikzlibrary{patterns}
\definecolor{dblue}{HTML}{1E90FF}
\definecolor{dred}{HTML}{bf2a2a}
\definecolor{black}{HTML}{bf2a2a}

\newtheoremstyle{define}{10pt}{0pt}{\itshape}{}{\bf}{.}{.5em}{}
\newtheoremstyle{exmp}{10pt}{0pt}{\itshape}{}{\bf}{.}{.5em}{}
\newtheoremstyle{rmrk}{10pt}{10pt}{}{}{\bf}{.}{.5em}{}

\newtheorem{theorem}{Theorem}[section]
\newtheorem{proposition}[theorem]{Proposition}

\newtheorem{lemma}[theorem]{Lemma}
\newtheorem{definition}[theorem]{Definition}
\theoremstyle{rmrk}
\newtheorem{remark}[theorem]{Remark}

\newcommand{\prf}[1]{\begin{proof}#1\end{proof}}
\newenvironment{hproof}{%
  \proof}{\endproof}

\title{Robust Persistence Diagrams \\using Reproducing Kernels}

\author{%
Siddharth Vishwanath \vspace{-1mm} \\  The Pennsylvania State University \vspace{-1mm} \\ {\url{suv87@psu.edu}}%
\And
Kenji Fukumizu$^*$ \vspace{-1mm} \\ The Institute of Statistical Mathematics \vspace{-1mm} \\  {\url{fukumizu@ism.ac.jp}}%
\And
Satoshi Kuriki$^*$ \vspace{-1mm} \\  The Institute of Statistical Mathematics \vspace{-1mm} \\  {\url{kuriki@ism.ac.jp}}%
\And
Bharath Sriperumbudur\thanks{Authors arranged alphabetically}\vspace{-1mm} \\  The Pennsylvania State University \vspace{-1mm} \\  {\url{bks18@psu.edu}}%
}

\begin{document}

\maketitle

\begin{abstract}
Persistent homology has become an important tool for extracting geometric and topological features from data, whose multi-scale features are summarized in a persistence diagram. From a statistical perspective, however, persistence diagrams are very sensitive to perturbations in the input space. In this work, we develop a framework for constructing robust persistence diagrams from superlevel filtrations of robust density estimators constructed using reproducing kernels. Using an analogue of the influence function on the space of persistence diagrams, we establish the proposed framework to be less sensitive to outliers. The robust persistence diagrams are shown to be consistent estimators in bottleneck distance, with the convergence rate controlled by the smoothness of the kernel---this in turn allows us to construct uniform confidence bands in the space of persistence diagrams. Finally, we demonstrate the superiority of the proposed approach on benchmark datasets.

\end{abstract}

\section{Introduction}
\label{introduction}

Given a set of points $\Xn = \pb{\Xv_1,\Xv_2,\dots,\Xv_n}$ observed from a probability distribution $\pr$ on an input space $\X \subseteq \R^d$, understanding the shape of $\Xn$ sheds important insights on low-dimensional geometric and topological features which underlie $\pr$, and this question has received increasing attention in the past few decades. To this end, Topological Data Analysis (TDA), with a special emphasis on persistent homology \citep{edelsbrunner2000topological,zomorodian2005computing}, has become a mainstay for extracting the shape information from data. In statistics and machine-learning, persistent homology has facilitated the development of novel methodology (e.g., \cite{chazal2013persistence,chen2019topological,bruel2018topology}), which has been widely used in a variety of applications dealing with massive, unconventional forms of data (e.g., \cite{bendich2016persistent,gameiro2015topological,xu2019finding}).

Informally speaking, persistent homology detects the presence of topological features across a range of resolutions by examining a nested sequence of spaces, typically referred to as a \textit{filtration}. The filtration encodes the birth and death of topological features as the resolution varies, and is presented in the form of a concise representation---a persistence diagram or barcode. In the context of data-analysis, there are two different methods for obtaining filtrations. The first is computed from the pairwise Euclidean distances of $\Xn$, such as the Vietoris-Rips, \cech{}, and Alpha filtrations \cite{edelsbrunner2000topological}. The second approach is based on 
choosing a function on $\X$ that reflects the density of $\pr$ (or its approximation based on $\Xn$), and, then, constructing a filtration. While the two approaches explore the topological features governing $\pr$ in different ways, in essence, they generate {similar} insights.

Despite obvious advantages, the adoption of persistent homology in mainstream statistical methodology is still limited. An important limitation among others, in the statistical context, is that the resulting persistent homology is highly sensitive to outliers. While the stability results of \cite{chazal2016structure,cohen2007stability} guarantee that small perturbations on all of $\Xn$ induce only small changes in the resulting persistence diagrams, a more pathological issue arises when a small fraction of $\Xn$ is subject to very large perturbations. Figure~\ref{fig:distfct} illustrates how inference from persistence diagrams can change dramatically when $\Xn$ is contaminated with only a few outliers.  Another challenge is the mathematical difficulty in performing sensitivity analysis in a formal statistical context. Since the space of persistence diagrams has an unusual mathematical structure, it falls victim to issues such as non-uniqueness of Fr\'{e}chet means and {unbounded} curvature of geodesics \citep{mileyko2011probability,turner2014frechet,divol2019understanding}. With this background, the \textit{central objective} of this paper is to develop outlier robust persistence diagrams, develop a framework for examining the sensitivity of the resulting persistence diagrams to noise, and establish statistical convergence guarantees. To the best of our knowledge, not much work has been carried out in this direction. {\citet{bendich2011improving} construct persistence diagrams from Rips filtrations on $\Xn$ by replacing the Euclidean distance with diffusion distance, \citet{brecheteau2018k} use a coreset of $\Xn$ for computing persistence diagrams from the distance-to-measure, and \citet{anai2019dtm} use weighted-Rips filtrations on $\Xn$ to construct more stable persistent diagrams.} However, no sensitivity analysis of the resultant diagrams are carried out in \citep{bendich2011improving,brecheteau2018k,anai2019dtm} to demonstrate their robustness.
\textbf{Contributions.} The main contributions of this work are threefold. 1) We propose robust persistence diagrams constructed from filtrations induced by an RKHS-based robust KDE (kernel density estimator) \cite{kim2012robust} of the underlying density function of $\pr$ (Section~\ref{robust-dgm}). While this idea of inducing filtrations by an appropriate function---\citep{fasy2014confidence,chazal2017robust,phillips2015geometric} use KDE, distance-to-measure (DTM) and kernel distance (KDist), respectively---has already been explored, we show the corresponding persistence diagrams to be less robust compared to our proposal.
2)
In Section~\ref{robustness}, we generalize the notions of \textit{influence function} and \textit{gross error sensitivity}---which are usually defined for normed spaces---to the space of persistence diagrams, which lack the vector space structure. Using these generalized notions, we investigate the sensitivity of persistence diagrams constructed from filtrations induced by different functions (e.g., KDE, robust KDE, DTM) and demonstrate the robustness of the proposed method, both mathematically (Remark~\ref{remark:influence}) and numerically (Section~\ref{experiments}). 
3) We establish the statistical consistency of the proposed robust persistence diagrams and provide uniform confidence bands by deriving exponential concentration bounds for the uniform deviation of the robust KDE (Section~\ref{consistency}).


\begingroup
\setlength{\intextsep}{-10pt} 
\setlength{\textfloatsep}{-10pt} 
\setlength{\abovecaptionskip}{-10pt}
\setlength{\belowcaptionskip}{-10pt}%
\begin{figure}[!t]
  \centering
  \begin{subfigure}[b]{0.24\linewidth}
    \includegraphics[width=\linewidth]{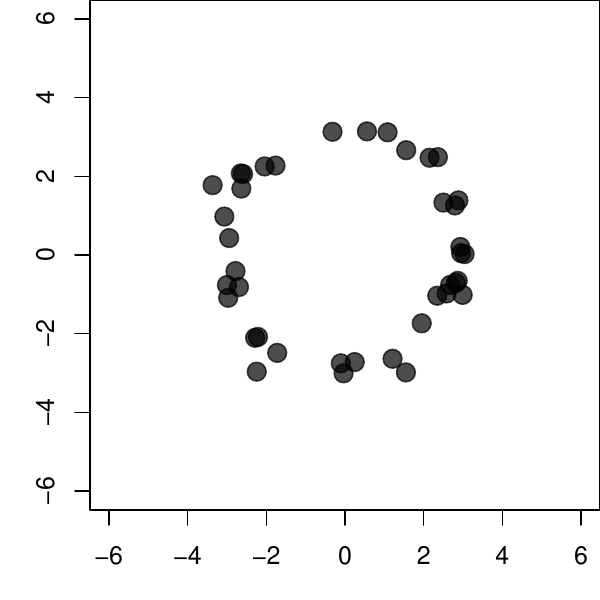}
  \end{subfigure}
  \begin{subfigure}[b]{0.24\linewidth}
    \includegraphics[width=\linewidth]{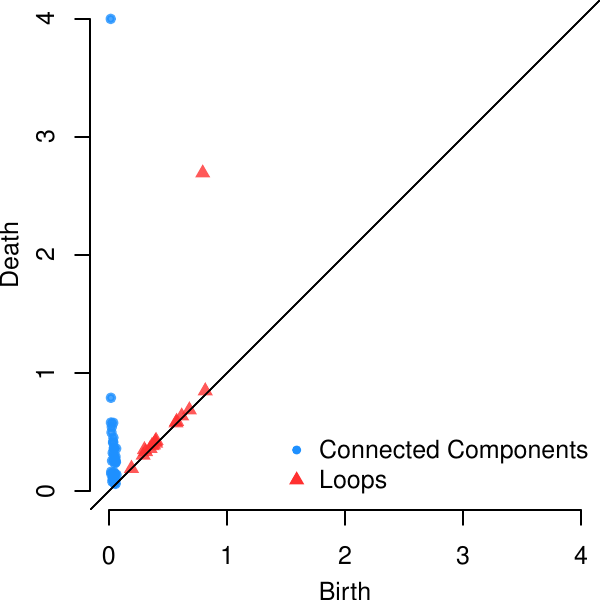}
  \end{subfigure}
  \begin{subfigure}[b]{0.24\linewidth}
    \includegraphics[width=\linewidth]{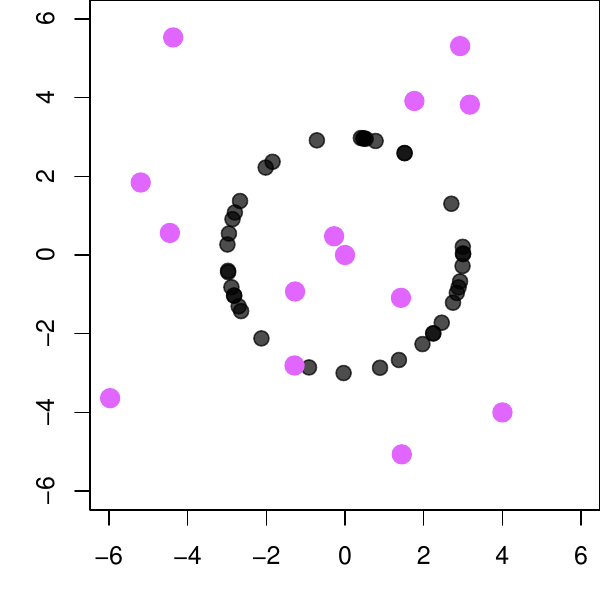}
  \end{subfigure}
  \begin{subfigure}[b]{0.24\linewidth}
    \includegraphics[width=\linewidth]{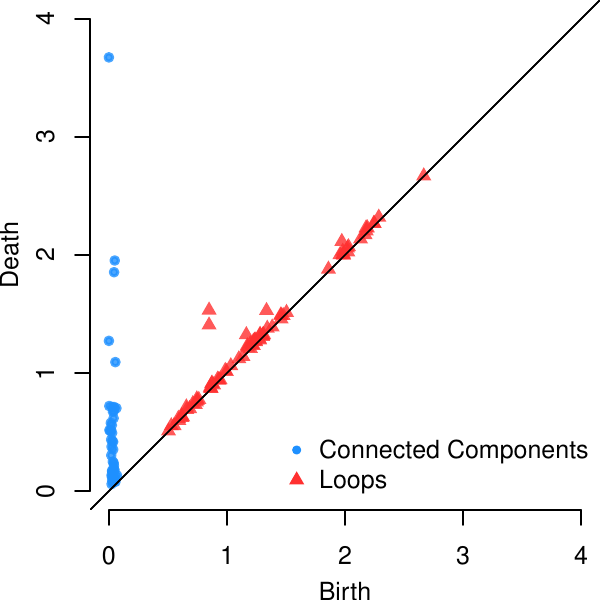}
  \end{subfigure}
  \caption{(Left) $\Xn$ is sampled from a circle with small perturbations to each point. The persistence diagram detects the presence of the loop, as guaranteed by the stability of persistence diagrams \cite{chazal2016structure,cohen2007stability}. (Right) $\Xn$ is sampled from a circle but with just a few outliers. The resulting persistence diagram changes dramatically --- the persistence of the main loop plummets, and other spurious loops appear, as elaborated in Section \ref{preliminaries}.}
  \label{fig:distfct}
\end{figure}
\endgroup

\textbf{Definitions and Notations.} {For a metric space $\X$, the ball of radius $r$ centered at $\xv\in \X$ is denoted by $B_{\X}(\xv,r)$. $\mathcal{P}(\R^d)$ is the set of all Borel probability measures on $\R^d$, and $\M(\R^d)$ denotes~the set of probability measures on $\R^d$ with compact support and \textit{tame} \mbox{density function (See Section \ref{preliminaries})}. $\delta_{\xv}$ denotes a Dirac measure \mbox{at $\xv$}. 
For bandwidth $\sigma > 0$, $\HH$ denotes a reproducing kernel Hilbert space (RKHS) with ${\K : \R^d\times\R^d \rightarrow \R}$ as its reproducing~kernel. We denote by $\Phis(\xv) = \K(\cdot,\xv) \in \HH$, the feature map associated with $\K$, which embeds $\xv \in \R^d$ into $\Phis(\xv) \in \HH$. Throughout this paper, we assume that $\K$ is radial, i.e., $\K (\xv,\yv) = \s^{-d}\psi({\norm{\xv-\yv}_2}/{\s})$ with $\psi(\Vert\cdot\Vert_2)$ being a pdf on $\R^d$, 
where $\Vert \xv\Vert^2_2=\sum^d_{i=1}x^2_i$ for ${\xv=(x_1,\ldots,x_d)\in\R^d}$. Some common examples include the Gaussian, Mat\'{e}rn and inverse multiquadric kernels. We denote ${\kinf \defeq \sup_{\xv,\yv\in \R^d}{\K(\xv,\yv)} = \s^{-d}\psi(0)}$. Without loss of generality, we assume $\psi(0)=1$. For $\pr \in \mathcal{P}(\R^d)$, $\mu_{\pr} \defeq \int\K(\cdot,\yv)d\pr(\yv) \in \HH$ is called the {mean embedding} of $\pr$, and ${\D \defeq \pb{\mu_{\pr} : \pr \!\in\!\mathcal{P}(\R^d)}}$ is the space of mean embeddings~\cite{muandet2017kernel}.}

\vspace{-2mm}
\section{Persistent Homology: Preliminaries}
\label{preliminaries}

We present the necessary background on persistent homology for completeness. See \cite{chazal2017introduction,wasserman2018topological} for a comprehensive introduction. 


\textbf{Persistent Homology.} Let $\phi : \X \rightarrow \R_{\geq 0}$ be a function on the metric space $\pa{\X,d}$. At level $r > 0$, the \textit{sublevel} set $\X_r = \phi\inv\pa{[0,r]} = \pb{\xv \in \X : \phi(\xv) \le r}$ encodes the \mbox{topological information} in $\X$. For $r < s$, the sublevel sets are nested, i.e., $\X_r \subseteq \X_s$. Thus $\pb{\X_r}_{0\le r < \infty}$ is a \mbox{nested sequence} of topological spaces, called a \textit{filtration}, denoted by $\text{Sub}(\phi)$, and $\phi$ is called the {\textit{filter function}}. As the level $r$ varies, the evolution of the topology is captured in the filtration. Roughly speaking, new~\mbox{cycles} (i.e.,~\mbox{connected components}, loops, voids and higher \mbox{order} analogues) can appear or existing cycles can merge. A new $k$-dimensional feature is said to be born at $b \in \R$ when a nontrivial \mbox{$k$-cycle} appears
\begin{minipage}[t]{0.5\textwidth}
    in $\X_b$. The same $k$-cycle dies at level ${d > b}$ when it disappears in all $\X_{d+\epsilon}$ for ${\epsilon > 0}$. \mbox{Persistent} \mbox{homology} is an algebraic module which tracks~the \textit{persistence pairs} $(b,d)$ of births $b$ and deaths~$d$ {with \mbox{multiplicity} $\mu$} across the entire filtration $\text{Sub}(\phi)$. \mbox{Mutatis mutandis}, a similar \mbox{notion} holds for {superlevel} sets ${\X^r = \phi\inv\pa{[r,\infty)}}$, inducing the {filtration} $\text{Sup}(\phi)$. For $r<s$, the inclusion ${\X^r \supseteq \X^s}$~is reversed and a cycle born at $b$ dies at a level $d<b$, resulting in the persistence pair $(d,b)$ instead. \mbox{Figure}~\ref{fig1} shows 3 connected {components} in the superlevel set for $r=8$. The components were born as $r$ swept through the blue points, and die when $r$ approaches the red points. In practice, \mbox{the filtrations are computed on a grid representation}
\end{minipage}\hfill
\begin{minipage}[t]{0.5\textwidth}
  \centering\raisebox{\dimexpr-\height+\ht\strutbox}{%
  \includegraphics[width=0.9\textwidth]{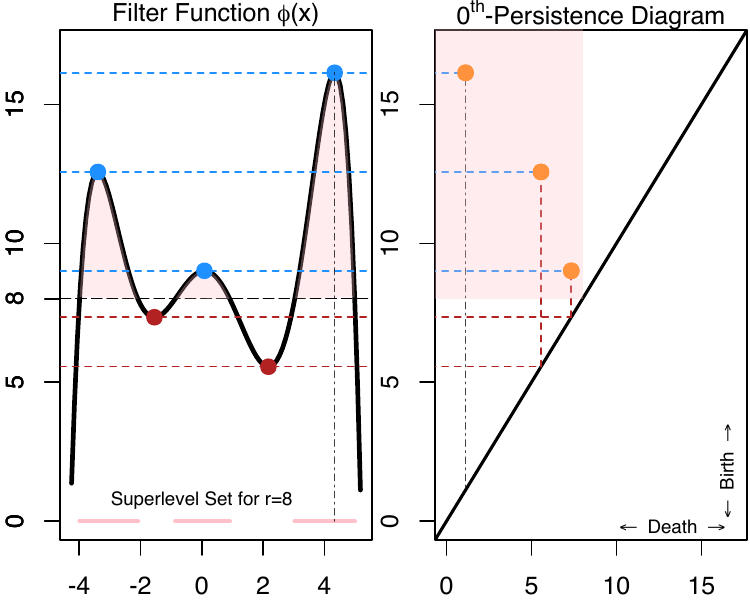}
  }
  \vspace{-1mm}
  \captionof{figure}{$\dgm\pa{\textup{Sup}(\phi)}$ for $\phi:\R\rightarrow \R$.}\vspace{1.5mm}
  \label{fig1}
\end{minipage}
  of the underlying space using cubical homology. We refer the reader to Appendix \ref{persistent-homology} for more details.

\textbf{Persistence Diagrams.} By collecting all persistence pairs, the persistent homology features are concisely represented as a persistence diagram $\dgm\pa{\text{Sub}(\phi)} \defeq \pb{(b,d) \in \R^2 : 0 \le b < d \le \infty}$. A similar definition carries over to $\dgm\pa{\text{Sup}(\phi)}$, using $(d,b)$ instead. See Figure~\ref{fig1} for an illustration. When the context is clear, we drop the reference to the filtration and simply write $\dgm(\phi)$. The $k^{th}$ persistence diagram is the subset of $\dgm(\phi)$ corresponding to the $k$-dimensional features. The space of persistence diagrams is the locally-finite multiset of points on ${\Omega = \pb{(x,y): 0 \le x < y \le \infty}}$, endowed with the family of $p$-Wasserstein metrics $W_p$, for $1 \le p \le \infty$. We \mbox{refer the} reader to \cite{edelsbrunner2010computational,divol2019understanding} for a thorough introduction. $W_\infty$ is commonly referred to as the \textit{bottleneck distance}.
\begin{definition}
  Given two persistence diagrams $D_1$ and $D_2$, the bottleneck distance is given by
  \eq{
  \Winf\pa{D_1,D_2} = \inf\limits_{\gamma \in \Gamma}\sup\limits_{p \in D_1 \cup \Delta} \norminf{p - \gamma(p)},\nn
  }
  where $\Gamma = \pb{\gamma : D_1 \cup \Delta \rightarrow D_2 \cup \Delta}$ is the set of all bijections from $D_1$ to $D_2$, including the diagonal $\Delta = \pb{(x,y) \in \R^2 : 0\le x=y \le \infty}$ with infinite multiplicity.\vspace{-2mm}
\end{definition}
An assumption we make at the outset is that the filter function $f$ is \textit{tame}. Tameness is a metric regularity condition which ensures that the number of points on the persistence diagrams are finite, and, in addition, the number of nontrivial cycles which share identical persistence pairings are also finite. Tame functions satisfy the celebrated stability property w.r.t. the bottleneck distance.
\begin{proposition}[Stability of Persistence Diagrams \cite{cohen2007stability,chazal2016structure}]
Given two tame functions ${f,g : \X \rightarrow \R}$,
$$\Winf\pa{\dgm(f),\dgm(g)} \le \norminf{f-g}.\nn$$
\textup{The space of persistence diagrams is, in general, challenging to work with. However, the stability property provides a handle on the persistence space through the function space of filter functions.
}
\label{prop:stability}
\end{proposition}
\vspace{-3mm}
\section{Robust Persistence Diagrams}
\label{robust-dgm}

Given $\Xn = \pb{\Xv_1,\Xv_2,\dots,\Xv_n} \subseteq \R^d$ drawn iid from a probability distribution $\pr \in \mathcal{M}(\R^d)$ with density $f$, the corresponding persistence diagram can be obtained by considering a \mbox{filter function} ${\phi_{n}:\R^d\rightarrow\R}$, constructed from $\Xn$ as an approximation to its population analogue, $\phi_\pr:\R^d\rightarrow\R$, that carries the topological information of $\pr$.

Commonly used $\phi_\pr$ include the (i) \mbox{kernelized density, $\barfo$,} (ii) Kernel Distance (KDist), $d^{\K}_{\pr}$, and (iii) distance-to-measure (DTM), $d_{\pr,m}$, which are defined as:
\eq{
\barfo(\xv) &\defeq \int_{\X}\K(\xv,\yv)d\pr(\yv) \ ;  & d^{\K}_{\pr} &\defeq \norm{\mu_{\delta_{\xv}} - \mu_{\pr}}_{\HH} \ ; & d_{\pr,m}(\xv) &\defeq \sqrt{\f{1}{m}\smallint_0^m{F\inv_{\xv}(u)du}},\nn
}
where $F_{\xv}(t) = \pr\pa{\norm{\Xv-\xv}_2 \le t}$ and $\sigma,m > 0$. For these $\phi_\pr$, the corresponding empirical analogues, $\phi_{n}$, are constructed by replacing $\pr$ with the empirical measure, $\pr_n\defeq\frac{1}{n}\sum^n_{i=1}\delta_{\Xv_i}$. For example, the empirical analogue of $\barfo$ is the familiar kernel density estimator (KDE), ${\barfn = \f{1}{n}\sum_{i=1}^n{\K(\cdot,\Xv_i)}}$.
While KDE and KDist encode {the shape and distribution of mass for $\supp({\pr})$} by approximating the density $f$ (sublevel sets of KDist are rescaled versions of superlevel sets of KDE \cite{phillips2015geometric,chazal2017robust}), DTM, on the other hand, approximates the distance function to $\supp({\pr})$.

Since $\phi_n$ is based on $\pr_n$, it is sensitive to outliers in $\Xn$, which, in turn affect the persistence diagrams (as illustrated in Figure \ref{fig:distfct}). To this end, in this paper, we propose \textit{robust persistence diagrams} constructed using superlevel filtrations of a robust density estimator of $f$, i.e., the filter function, $\phi_n$ is chosen to be a robust density estimator of $f$. Specifically, we use the robust KDE, $\fn$, introduced by \cite{kim2012robust} as the filter function, which is defined as a solution to the following
M-estimation problem:
\eq{
\fn \defeq \arginf\limits_{g\in \G} \int_{\X}{\rho\pa{\norm{\Phis(\yv) - g}_{\HH}} d\pr_n(\yv)},
\label{def:rkde}
}
where $\rho : \R_{\ge 0} \rightarrow \R_{\ge 0}$ is a robust loss function, and $\G = \HH \cap \D = \D$ is the hypothesis class. Observe that when $\rho(z) = \f{1}{2}z^2$, the unique solution to Eq.~\eqref{def:rkde} is given by the KDE, $\barfn$. 
Therefore, a robust KDE is obtained by replacing the square loss with a \textit{robust loss}, which satisfies the following assumptions. These assumptions, which are similar to those of \cite{kim2012robust,vandermeulen2013consistency} guarantee the existence and uniqueness (if $\rho$ is convex) of $\fn$ \cite{kim2012robust}, and are satisfied by most robust loss functions, including the Huber loss, $\rho(z) = \f{1}{2}z^2\mathbbm{1}\pb{z \le 1} + \pa{z-\half}\mathbbm{1}\pb{z>1}$ and the Charbonnier loss, $\rho(z) = \sqrt{1+z^2}-1$.
\vspace{-.5mm}
\begin{description}
  \item[$\pa{\mathcal{A}1}$] $\rho$ is strictly-increasing and $M$-Lipschitz, with $\rho(0)=0$.
  \item[$\pa{\mathcal{A}2}$] $\rho'(x)$ is continuous and bounded with $\rho'(0)=0$ .
  \item[$\pa{\mathcal{A}3}$] $\varphi(x) = \rho'(x)/x$ is bounded, $L$-Lipschitz and continuous, with $\varphi(0) < \infty$.
  \item[$\pa{\mathcal{A}4}$] $\rho''$ exists, with $\rho''$ and $ \varphi$ nonincreasing.\vspace{-.5mm}
\end{description}
Unlike for squared loss, the solution $\fn$ cannot be obtained in a closed form. However, it can be shown to be the fixed point of an iterative procedure, referred to as KIRWLS algorithm \cite{kim2012robust}. The KIRWLS algorithm starts with initial weights $\{w^{(0)}_i\}^{n}_{i=1}$ such that $\sum_{i=1}^{n}w_i^{(0)}=1$, and generates the iterative sequence of estimators $\{f_{\rho,\s}^{(k)}\}_{k\in\mathbb{N}}$ as
\eq{
f_{\rho,\s}^{(k)} = \sum_{i=1}^n{w_i^{(k-1)}\K(\cdot,\Xv_i)} \ \ ; && w_i^{(k)}= \f{\varphi(\Vert\Phis(\Xv_i)-f_{\rho,\s}^{(k)}\Vert_{\HH})}{\sum_{j=1}^{n}{\varphi(\Vert\Phis(\Xv_j)-f_{\rho,\s}^{(k)}\Vert_{\HH})}}.&\nn
}
Intuitively, note that if $\Xv_i$ is an outlier, then the corresponding weight $w_i$ is small (since $\varphi$ is nonincreasing) and therefore less {weight} is given to the contribution of $\Xv_i$ in the density estimator. Hence, the weights serve as a measure of \emph{inlyingness}---smaller (\textit{resp.} larger) the weights, lesser (\textit{resp.} more) inlying are the points.
When $\pr_n$ is replaced by $\pr$, the solution of Eq.~\eqref{def:rkde} is its population analogue, $\fo$. Although $\fo$ does not admit a closed form solution, it can be shown \cite{kim2012robust} that there exists a non-negative real-valued function $\ws$ satisfying $\int_{\R^d} \ws(\xv)\,d\pr(\xv)=1$ such that
\eq{
\fo = \int_{\R^d}\K(\cdot,\xv)\ws(\xv)d\pr(\xv) = \int_{\R^d}{\f{\varphi ({\normh{\phifox}}) }{\int_{\R^d}\varphi ({\normh{\Phis(\yv)-\fo}}) d\pr(\yv)} \K(\cdot,\xv) \ d\pr(\xv) },
\label{eq:weight}
}
where $\ws$ acts as a population analogue of the weights in KIRWLS algorithm.

To summarize our proposal, the fixed point of the KIRWLS algorithm, which yields the robust density estimator $\fn$, is used as the filter function to obtain a robust persistence diagram of $\Xn$. On the computational front, note that $\fn$ is computationally more complex than the KDE, $\barfn$, requiring $O(n\ell)$ computations compared to $O(n)$ of the latter, with $\ell$ being the number of iterations required to reach the fixed point of KIRWLS. However, once these filter functions are computed, the corresponding persistence diagrams have similar computational complexity as both require computing superlevel sets, which, in turn, require function evaluations that scale as $O(n)$ for both $\fn$ and $\barfn$.

\section{Theoretical Analysis of Robust Persistence Diagrams}
\label{results}

In this section, we investigate the theoretical properties of the proposed robust persistence diagrams. First, in Section~\ref{robustness}, we examine the sensitivity of persistence diagrams to outlying perturbations through the notion of \textit{metric derivative} and compare the effect of different filter functions.
Next, in Section~\ref{consistency}, we establish consistency and convergence rates for the robust persistence diagram to its population analogue. These results allow to construct uniform confidence bands for the robust persistence diagram. The proofs of the results are provided in Appendix~\ref{proofs}.


\subsection{A measure of sensitivity of persistence diagrams to outliers}
\label{robustness}

The influence function and gross error sensitivity are arguably the most popular tools in robust statistics for diagnosing the sensitivity of an estimator to a single adversarial contamination~\cite{hampel2011robust,huber2004robust}. Given a statistical functional $T:\mathcal{P}(\X) \rightarrow \pa{V,\norm{\cdot}_V}$, which takes an input probability measure $\pr \in \mathcal{P}(\X)$ on the input space $\X$ and produces a statistic $\pr \mapsto T(\pr)$ in some normed space $\pa{V,\norm{\cdot}_V}$, the \textit{influence function} of $\xv \in \X$ at $\pr$ is given by the \gat{} derivative of $T$ at $\pr$ restricted to the space of signed Borel measures with zero expectation:\vspace{-1.5mm}
\eq{
\mathsf{IF}(T;\pr,\xv) \defeq \f{\partial}{\partial \epsilon}T \Bigl( (1-\epsilon)\pr + \epsilon \delta_{\xv} \Bigr) \Big\vert_{\epsilon=0} = \lim_{\epsilon \rightarrow 0} \f{T\pa{(1-\epsilon)\pr + \epsilon \delta_{\xv}} - T(\pr)}{\epsilon},\nn
}
and the \textit{gross error sensitivity} at $\pr$ is given by $\Gamma(T;\pr) \defeq \sup_{\xv \in \X} \norm{\mathsf{IF}(T;\pr,\xv)}_V$. However, a persistence diagram (which is a statistical functional) does not take values in a normed space and therefore the {notion of influence functions} has to be generalized to metric spaces through the concept of a metric derivative: Given a complete metric space $(X,d_X)$ and a curve ${s : [0,1] \rightarrow X}$, the \textit{metric derivative} at $\epsilon=0$ is given by
$\abs{s'}(0) \defeq \lim_{\epsilon \rightarrow 0}\frac{1}{\epsilon}d_X(s(0),s(\epsilon)).$
Using this generalization, we have the following definition, which allows to examine the influence an outlier has on the persistence diagram obtained from a filtration.
\begin{definition}
  Given a probability measure $\pr \in \mathcal{P}(\R^d)$ and a filter function $\phi_{\pr}$ depending on $\pr$, the {\em persistence influence} of a perturbation $\xv \in \R^d$ on $\dgm\pa{\phi_{\pr}}$ is defined as
  \eq{
  \Psi\pa{\phi_{\pr};\xv} = \lim_{\epsilon \rightarrow 0}\f{1}{\epsilon} \Winf\pa{\dgm\pa{\phi_{\pr^\epsilon_{\xv}}},\dgm\pa{\phi_\pr}},\nn
  }
  where $\pr^\epsilon_{\xv} \defeq (1-\epsilon)\pr + \epsilon \delta_{\xv}$, and the {\em gross-influence} is defined as ${\Gamma(\phi_{\pr}) = \sup_{\xv \in \R^d}\Psi\pa{\phi_{\pr};\xv}}$.\vspace{-.5mm}
\end{definition}%
{For $\e > 0$, let $\fox$ be the robust KDE associated with the probability measure $\pr^\epsilon_{\xv}$.} The following result (proved in Appendix~\ref{proof:influence}) bounds the persistence influence for the persistence diagram induced by the filter function $\fo$, which is the population analogue of robust KDE.
\begin{theorem}
  For a loss $\rho$ satisfying $\pa{\mathcal{A}1}$--$\pa{\mathcal{A}3}$, and $\sigma > 0$, if { $\lim\limits_{\epsilon \rightarrow 0}\f{1}{\epsilon}\pa{\fox-\fo}$ exists}, then the persistence influence of $\xv\in\R^d$ on $\dgm\pa{\fo}$ satisfies
  \eq{
  \Psi\pa{\fo;\xv} \le \kinf^\half  \rho'\pa{\normh{\phifox}}\left(\int_{\R^d}\zeta{\pa{\normh{\phifo}}}d\pr(\yv)\right)^{-1},
  \label{eq:influence-rkde}
  }
  where $\zeta(z) = \varphi(z) - z\varphi'(z)$.
  %
  \label{thm:influence}
\end{theorem}
\begin{remark}
We make the following observations from Theorem \ref{thm:influence}.

\quad \textbf{(i)} Choosing $\rho(z) = \f{1}{2}z^2$ and noting that $\varphi(z) = \rho''(z) = 1$, {a similar analysis, as in the proof of Theorem~\ref{thm:influence}, yields a bound} for the persistence influence of the KDE as
$$
\Psi\pa{\barfo;\xv} \le \sigma^{-d/2}
\normh{\Phis(\xv) - \barfo}.
$$
On the other hand, for robust loss functions, the term in Eq.~\eqref{eq:influence-rkde} involving $\rho'$ is bounded because of $\pa{\mathcal{A}2}$, making them less sensitive to very large perturbations. In fact, for nonincreasing $\varphi$, it can be shown (see Appendix~\ref{influence}) 
that%
\eq{
\Psi\pa{\fo;\xv} \le \sigma^{-d/2}\ws(\xv)
\normh{\Phis(\xv) - \fo},\nn
}
where, in contrast to KDE, %
the measure of inlyingness, $\ws$, weighs down extreme outliers.\vspace{-1mm}

\quad \textbf{(ii)} For the generalized Charbonnier loss (a robust loss function), given by ${\rho(z) = \pa{1+z^2}^{\alpha/2}-1}$ for $1\le \alpha < 2$, the persistence influence satisfies
\eq{
\Psi\pa{\fo;\xv} \le \sigma^{-d/2} \pa{1+\normh{\Phis(\xv) - \fo}^2}^{\f{\alpha-1}{2}} \pa{1+\int_{\R^d}\normh{\Phis(\yv) - \fo}^2 d\pr(\yv)}^{\f{1-\alpha}{2}}.\nn
}
 Note that for $\alpha=1$, the bound on the persistence influence $\Psi\pa{\fo;\xv}$ does not depend on how extreme the outlier $\xv$ is. Similarly, for the Cauchy loss, given by $\rho(z) = \log(1+z^2)$, we have
\eq{
\Psi\pa{\fo;\xv} \le \sigma^{-d/2} \pa{1+\int_{\R^d}\normh{\Phis(\yv) - \fo}^2d\pr(\yv)}.\nn
}
 This shows that for large perturbations, the gross error sensitivity for the Cauchy and \mbox{Charbonnier} losses are far more stable than that of KDE. This behavior is
also empirically illustrated in Figure \ref{fig:influence}. The experiment is detailed in Appendix~\ref{influence}.%

  \quad \textbf{(iii)} For the DTM function, it can be shown that   \eq{
  \Psi\pa{d_{\pr,m};\xv} \le \f{2}{\sqrt{m}}\sup\pb{\Abs{f(\xv) - \int_{\R^d}f(\yv)d\pr(\yv)} : \norm{\nabla f}_{L_2(\pr)} \le 1}.
  \label{eq:influence-dtm}
  }
  While $d_{\pr,m}$ cannot be compared to both $\barfo$ and $\fo$, as it captures topological information at a different scale, determined by $m$, we point out that when $\emph{supp}(\pr)$ is compact, $\Psi\pa{d_{\pr,m};\xv}$ is not guaranteed to be bounded, unlike in $\Psi\pa{\fo;\xv}$. We refer the reader to Appendix \ref{influence} for more details.\vspace{-1.5mm}
  \label{remark:influence}
\end{remark}

\begin{figure}[t]
 \centering
 \begin{subfigure}[b]{0.32\linewidth}
   \includegraphics[width=\linewidth]{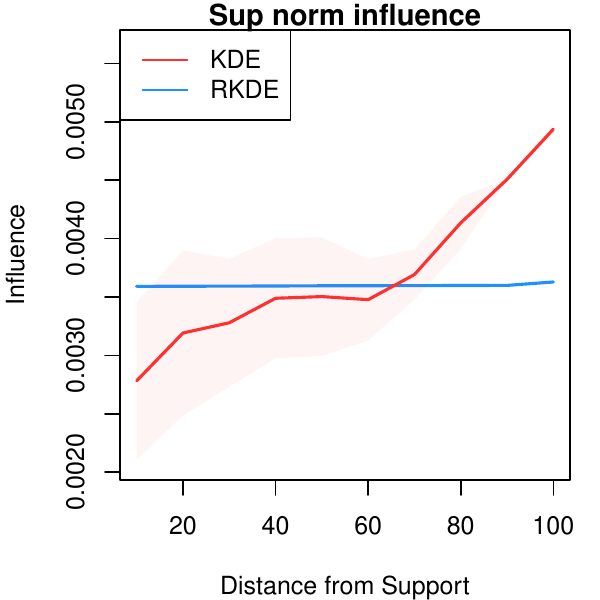}
 \end{subfigure}
 \begin{subfigure}[b]{0.32\linewidth}
   \includegraphics[width=\linewidth]{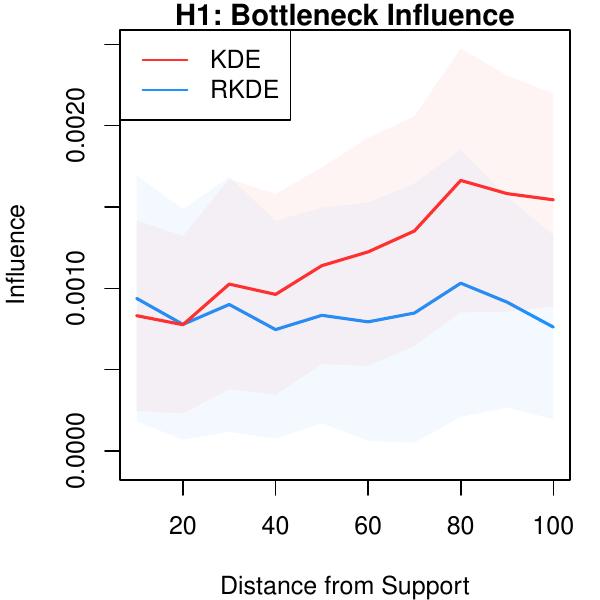}
 \end{subfigure}
 \begin{subfigure}[b]{0.32\linewidth}
   \includegraphics[width=\linewidth]{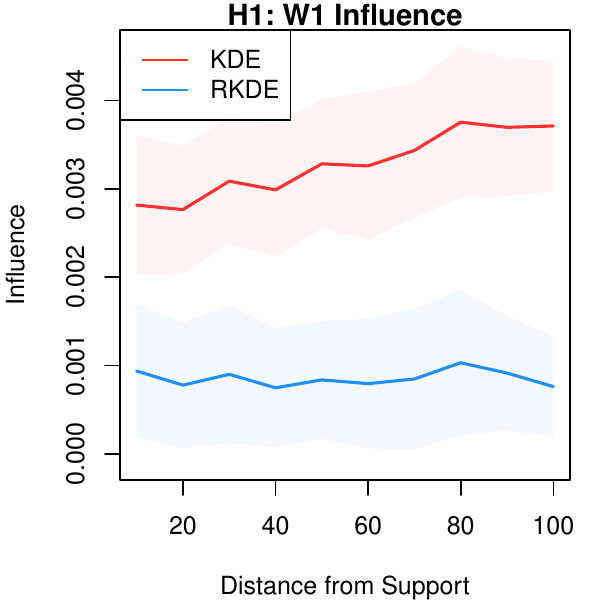}
 \end{subfigure}
 \vspace{-0.1mm}
 \caption{Points $\Xn$ are sampled from $\pr$ with nontrivial $1^{st}$-order homological features and outliers $\Y_m$ are added at a distance $r$ from the support of $\pr$. (Left) The average $L_{\infty}$ distance between the density estimators computed using $\Xn$ and $\Xn\cup\Y_m$ as $r$ increases. (Center) The average $W_{\infty}$ distance between the corresponding persistence diagrams for the $1^{st}$-order homological features. (Right) The $W_1$ distance (defined in \mbox{Eq.~\ref{eq:p-wasserstein} in Appendix~\ref{persistent-homology}}) between the same persistence diagrams. The results show that the outliers $\Y_m$ have little influence on the persistence diagrams from the robust KDEs. In contrast, as the outliers become more extreme (i.e., $r$ increases) their influence on the persistence diagrams from the KDE becomes more prominent.}
 \label{fig:influence}
 \end{figure}

It follows from Remark~\ref{remark:influence} that as $\s \rightarrow 0$,
the persistence influence of both the KDE and robust KDE behave as $O(\sigma^{-d})$,
showing that the robustness of robust persistence diagrams manifests only in cases where $\s > 0$. However, robustness alone has no bearing if the robust persistence diagram and the persistence diagram from the KDE are fundamentally different, i.e., they estimate different quantities as $\s \rightarrow 0$. The following result (proved in Appendix~\ref{proof:consistency2}) shows that as $\sigma \rightarrow 0$, $\dgm\pa{\fo}$ recovers the same information as that in  $\dgm\pa{\barfo}$, which is same as $\dgm\pa{f}$, where $f$ is the \mbox{density of $\pr$}.
\begin{theorem}
  For a strictly-convex loss $\rho$ satisfying $\pa{\mathcal{A}1}$--$\pa{\mathcal{A}4}$, and $\s > 0$, suppose $\pr \in \M(\R^d)$ with density $f$, and $\fo$ is the robust KDE. Then $\Winf\pa{\dgm\pa{\fo},\dgm\pa{f}} \rightarrow 0$ as $\sigma \rightarrow 0$.\vspace{-.5mm}
  \label{thm:consistency2}
\end{theorem}

Suppose $\pr = (1-\pi)\pr_0 + \pi \qr$, where $\pr_0$ corresponds to the true signal which we are interested in studying, and $\qr$ manifests as some ambient noise with $0 < \pi < \half$. In light of Theorem \ref{thm:consistency2}, by letting $\sigma \rightarrow 0$, along with the topological features of $\pr_0$, we are also capturing the topological features \mbox{of $\qr$}, which may obfuscate any statistical inference made using the persistence diagrams. In a manner, choosing $\sigma > 0$ suppresses
the noise 
in the resulting persistence diagrams, thereby making them more stable. On a similar note, the authors in \cite{fasy2014confidence} state that for a suitable bandwidth $\sigma > 0$, the level sets of $\barfo$ carry the same topological information as $\supp(\pr)$, despite the fact that some subtle details in $f$ may be omitted. In what follows, we consider the setting where robust persistence diagrams are constructed for a fixed $\sigma > 0$.
%
%


\subsection{Statistical properties of robust persistence diagrams from samples}
\label{consistency}
Suppose $\dgm\pa{\fn}$ is the robust persistence diagram obtained from the robust KDE on a sample $\Xn$ and $\dgm\pa{\fo}$ is its population analogue obtained from $\fo$.  
The following result (proved in Appendix~\ref{proof:consistency}) establishes the consistency of $\dgm\pa{\fn}$ in the~$\Winf$~metric.
\smallskip
\begin{theorem}
  For convex loss $\rho$ satisfying $\pa{\mathcal{A}1}$--$\pa{\mathcal{A}4}$, and fixed $\sigma>0$, suppose $\Xn$ is observed iid from a distribution $\pr\!\in\!\mathcal{M}(\R^d)$ with density $f$. Then
  \eq{
  {{W}_{\infty}\pa{\dgm\pa{\fn},\dgm\pa{\fo}} \stackrel{p}{\rightarrow} 0\  \text{ \ \ as } n \rightarrow \infty.}\nn
  }
  \vspace{-3mm}
  \label{thm:consistency}
\end{theorem}
We present the convergence rate of the above convergence in Theorem~\ref{thm:concentration1}, which depends on the smoothness of $\HH$. In a similar spirit to \cite{fasy2014confidence}, this result paves the way for constructing uniform confidence bands. Before we present the result, we first introduce the notion of \textit{entropy numbers} associated with an RKHS.


\begin{definition}[Entropy Number]
  Given a metric space $\pa{T,d}$ the $n^{th}$ entropy number is defined as
  \eq{
  e_n(T,d) \defeq \inf\pb{\epsilon > 0 : \exists \ \pb{t_1, t_2, \dots , t_{2^{n-1}}} \subset T  \ \text{ such that } T \subset \mathop{\bigcup}_{\scriptsize i=1}^{\scriptsize 2^{n-1}}B_{d}(t_i,\epsilon) }.\nn
  }
  Further, if $\pa{V,\norm{\cdot}_V}$ and $\pa{W,\norm{\cdot}_W}$ are two normed spaces and $L: V \rightarrow W$ is a bounded, linear operator, then $e_n(L) = e_n(L: V \rightarrow W) \defeq e_n\pa{L(B_V),\norm{\cdot}_W}$, where $B_V$ is a unit ball in $V$.
\end{definition}
Loosely speaking, entropy numbers are related to the eigenvalues of the integral operator associated with the kernel $\K$, and measure the capacity of the RKHS in approximating functions in $L_2(\R^d)$. In our context, the entropy numbers will provide useful bounds on the covering numbers of sets in the hypothesis class $\G$. We refer the reader to \cite{steinwart2008support} for more details. With this background, the following theorem (proved in Appendix~\ref{proof:concentration}) provides a method for constructing uniform confidence bands for the persistence diagram constructed using the robust KDE on $\Xn$.
\begin{theorem}
  For convex loss $\rho$ satisfying $\pa{\mathcal{A}1}$--$\pa{\mathcal{A}4}$, and fixed $\sigma > 0$, suppose the kernel $\K$ satisfies ${e_n\pa{\id : \HH \rightarrow L_{\infty}(\X)} \le a_\s n^{-\f{1}{2p}}}$, where $a_\s>1$, $0<p<1$ and $\X\subset\R^d$. Then, for a fixed \mbox{confidence level $0 < \alpha < 1$},
  \eq{
  \sup_{\pr \in \mathcal{M}(\X)}\pr^{\otimes n}\Bigg\{\Winf \Big( \dgm\left(\fn\right),\dgm\left(\fo\right) \Big) > \f{2M \kinf^\half}{\mu}\pa{\xi(n,p) + \delta\sqrt{\f{2\log\pa{{1}/{\alpha}}}{n}}} \Bigg\} \le \alpha,\nn
  }
  where $\xi(n,p)$ is given by%
  \eq{
  \xi(n,p) = \begin{cases}\medskip
  \gamma \f{a_\s^p}{(1-2p)} \cdot \f{1}{\sqrt{n}} & \textup{ if } 0 < p < 1/2,\nn\\\medskip
  {\gamma C\sqrt{a_\s}}\cdot\f{\log(n)}{\sqrt{n}} & \textup{ if }  p = 1/2,\nn\\\medskip 
  \gamma\f{p\sqrt{a_\s}}{2p-1} \cdot \f{1}{n^{{1}/{4p}}} & \textup{ if } 1/2 < p < 1,\nn
\end{cases}
}
for fixed constants  ${\gamma > \frac{12}{\sqrt{\log{2}}}, \ {C > 3 - \log(9a_\s)}}$ and $\mu = 2\min\pb{\varphi(2\kinf^\half), \rho''{(2\kinf^\half)}}$.
\label{thm:concentration1}
\end{theorem}


\begin{remark}

  We highlight some salient observations from Theorem \ref{thm:concentration1}.

  (i) If
  $diam(\X) = r$, and the kernel $\K$ is $m$-times differentiable, then from \cite[Theorem 6.26]{steinwart2008support}, the entropy numbers associated with $\K$ satisfy ${e_n\pa{\id : \HH \rightarrow L_{\infty}(\X)} \le c r^m n^{-\f{m}{d}}}$. In light of Theorem \ref{thm:concentration1}, for $p = \f{d}{2m}$, we can make two important observations.  First, as the dimension of the input space $\X$ increases, we have that the rate of convergence decreases; which is a direct consequence from the curse of dimensionality. Second, for a fixed dimension of the input space, the parameter $p$ in Theorem \ref{thm:concentration1} can be understood to be inversely proportional to the smoothness of the kernel. Specifically, as the smoothness of the kernel increases, the rate of convergence is faster, and we obtain sharper confidence bands. This makes a case for employing smoother kernels.\vspace{-.5mm}

  (ii) A similar result is obtained in \cite[Lemma 8]{fasy2014confidence} for persistence diagrams from the KDE, with a convergence rate $O_p({n^{-1/2}})$, where the proof relies on a simple application of Hoeffding's inequality, unlike the sophisticated tools the proof of Theorem \ref{thm:concentration1} warrants for the robust KDE.
\end{remark}



\setcounter{footnote}{0}

\vspace{-5mm}
\section{Experiments}
\label{experiments}

We illustrate the performance of robust persistence diagrams in machine learning applications through synthetic and real-world experiments.\footnote{\url{https://github.com/sidv23/robust-PDs}} In all the experiments, the kernel bandwidth $\s$ is chosen as the median distance of each $\xv_i \in \Xn$ to its $k^{th}$--nearest neighbour using the Gaussian kernel with the Hampel loss (similar setting as in \cite{kim2012robust})---we denote this bandwidth as $\s(k)$. Since DTM is closely related to the $k$-NN density estimator \cite{biau2011weighted}, we choose the DTM smoothing parameter as $m(k)={k}/{n}$. Additionally, the KIRWLS algorithm is run until the relative change of empirical risk $<10^{-6}$.

{\textbf{Runtime Analysis.} For $n=1000$, $\Xn$ is sampled from a torus inside $[0,2]^3$. For each grid resolution $\alpha \in \pb{0.04, 0.06, 0.08, 0.10}$, the robust persistence diagram $\dgm\pa{\fn}$ and the KDE persistence diagram $\dgm\pa{\fns}$ are constructed from the superlevel filtration of cubical homology. The \textit{total} time taken to compute the persistence diagrams is reported in Table \ref{tab:runtime}. The results demonstrate that the computational bottleneck is the persistent homology pipeline, and not the KIRWLS for $\fn$.}
\vspace{-3mm}

\begin{table}[H]
\caption{Runtime (in Seconds) for computing $\dgm\pa{\fn}$ and $\dgm\pa{\fns}$ at each grid resolution.\vspace{-2mm}
\label{tab:runtime}}
\small
\begin{tabu}{|l|c|c|c|c|}\hline
  Grid Resolution         & $0.04$ & $0.06$ & $0.08$ & $0.10$ \\ \tabucline[2pt]{-}
  Average runtime for $\dgm\pa{\fn}$             & $76.7$s & $17.1$s & $6.7$s & $3.5$s \\ \hline
  Average runtime for $\dgm\pa{\fns}$             & $75.5$s & $15.3$s & $4.7$s & $1.8$s \\ \hline
\end{tabu}
\end{table}
\vspace{-4mm}

\begin{figure}
  \centering
  \begin{subfigure}[b]{0.24\linewidth}
    \includegraphics[width=\linewidth]{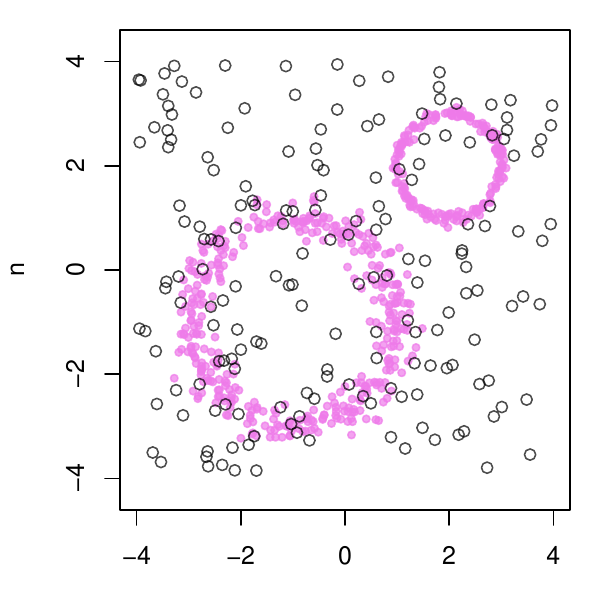}
    \caption{$\Xn$ (in \textcolor{magenta}{$\bullet$}) and $\Y_m$ (in $\circ$)}
  \end{subfigure}
  \begin{subfigure}[b]{0.24\linewidth}
    \includegraphics[width=\linewidth]{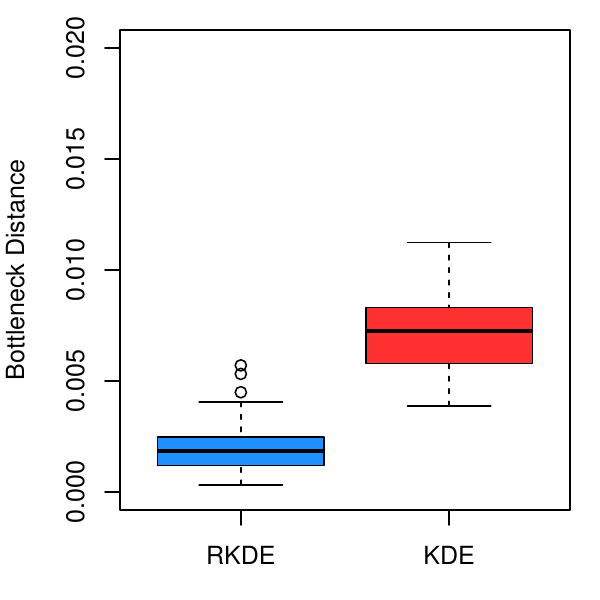}
    \caption{$\pi=20\%$, $p=4\times10^{-60}$}
  \end{subfigure}
  \begin{subfigure}[b]{0.24\linewidth}
    \includegraphics[width=\linewidth]{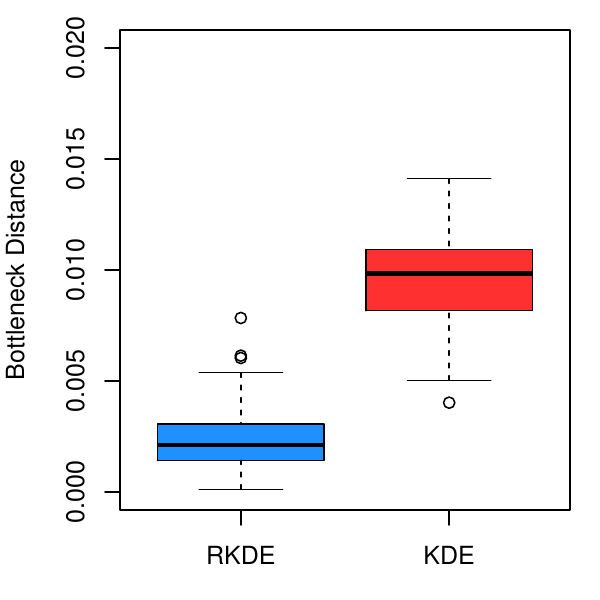}
    \caption{$\pi=30\%$, $p=2\times10^{-72}$}
  \end{subfigure}
  \begin{subfigure}[b]{0.24\linewidth}
    \includegraphics[width=\linewidth]{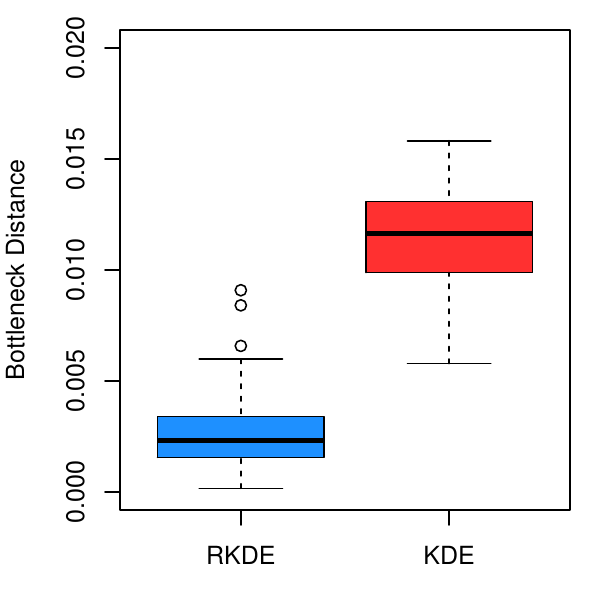}
    \caption{$\pi=40\%$, $p=2.5\times10^{-75}$}
  \end{subfigure}
   \vspace{-0.5mm}
  \caption{(a) A realization of $\Xn \cup \Y_m$. (b, c, d) As the noise level $\pi$ increases, boxplots for $\Winf\pa{\mathbf{D}_{\rho,\s},\mathcal{D}_\s^\#}$ in blue and $\Winf\pa{\mathbf{D}_{\s},\mathcal{D}_\s^\#}$ in red show that the robust persistence diagram recovers the underlying signal better.}
  \label{fig:bottleneck}
  \vspace{-0.5mm}
\end{figure}

\textbf{Bottleneck Simulation.} The objective of this experiment is to assess how the robust KDE persistence diagram compares to the KDE persistence diagram in recovering the topological features of the underlying signal. $\Xn$ is observed uniformly from two circles and $\Y_m$ is sampled uniformly from the enclosing square such that $m=200$ and ${m}/{n} = \pi \in \pb{20\%,30\%,40\%}$---shown in Figure~\ref{fig:bottleneck} (a). For each noise level $\pi$, and for each of $N=100$ realizations of $\Xn$ and $\Y_m$, the robust persistence diagram $\mathbf{D}_{\rho,\sigma}$ and the KDE persistence diagram $\mathbf{D}_\s$ are constructed from the noisy samples $\Xn \cup \Y_m$. In addition, we compute the KDE persistence diagram $\mathcal{D}_\s^\#$ on $\Xn$ alone as a proxy for the target persistence diagram one would obtain in the absence of any contamination. The bandwidth $\s(k) > 0$ is chosen for $k=5$. For each realization $i$, bottleneck distances ${U_i= \Winf\pa{\mathbf{D}_{\rho,\s},\mathcal{D}_\s^\#}}$ and ${V_i=\Winf\pa{\mathbf{D}_{\s},\mathcal{D}_\s^\#}}$ are computed for $1^{st}$-order homological features. The boxplots and $p$-values for the one-sided hypothesis test $H_0: U-V = 0$ vs. $H_1:U-V<0$ are reported in Figures~\ref{fig:bottleneck}~(b,~c,~d). The results demonstrate that the robust persistence diagram is noticeably better in recovering the true homological features, and in fact demonstrates superior performance when the noise levels are higher.

{\textbf{Spectral Clustering using Persistent Homology.} We perform a variant of the six-class benchmark experiment from \cite[Section 6.1]{adams2017persistence}. The data comprises of six different $3$D ``objects'': \texttt{cube, circle, sphere, 3clusters, 3clustersIn3clusters,} and \texttt{torus}. $25$ point clouds are sampled from each object with additive Gaussian noise (SD$= 0.1$), \textit{and} ambient Mat\'ern cluster noise. For each point cloud, $\Xn$, the robust persistence diagram $\dgm\pa{\fn}$ and the persistence diagram $\dgm\pa{d_{\Xn}}$, from the {distance function}, are constructed. Additionally, $\dgm\pa{d_{\Xn}}$ is transformed to the persistence image $\Img\pa{d_{\Xn}, h}$ for $h=0.1$. Note that $\dgm\pa{\fn}$ is a robust diagram while $\Img\pa{d_{\Xn}, h}$ is a stable vectorization of a non-robust diagram \cite{adams2017persistence}. For each homological order $\pb{H_0, H_1, H_2}$, distance matrices $\pb{\Delta_0, \Delta_1, \Delta_2}$ are computed: $W_p$ metric for $\dgm\pa{\fo}$, and $L_p$ metric for $\Img\pa{d_{\Xn},h}$ with $p \in \pb{1, 2, \infty}$, and spectral clustering is performed on the resulting distance-matrices. The quality of the clustering is assessed using the rand-index. The results, reported in Table \ref{tab:spectral}, evidence the superiority of employing inherently robust persistence diagrams in contrast to a robust vectorization of an inherently noisy persistence diagram.}
\vspace{-2mm}

\begin{table}[H]
\caption{Rand-index for spectral clustering using distance matrices for $\dgm\pa{\fo}$ and $\Img\pa{d_{\Xn},h}$.\vspace{-2mm}
\label{tab:spectral}}
\small
{\renewcommand{\arraystretch}{1.2}
\begin{tabu}{|l|c c c|c c c|} \cline{2-7}
  \multicolumn{1}{c|}{} &  \multicolumn3{c|}{$\dgm\pa{\fo}$} & \multicolumn3{c|}{$\Img\pa{d_{\Xn},h}$} \\ \tabucline-
  \tabuphantomline
  Distance Metric & $W_1$ & $W_2$ & $\Winf$ & $L_1$ & $L_2$ & $L_{\infty}$ \\ \tabucline[2pt]{-}\tabuphantomline
$\Delta_0$ (from $H_0$)  & $95.30\%$ & $93.65\%$ & $94.44\%$ & $78.53\%$ & $81.77\%$ & $80.05\%$ \\
$\Delta_1$ (from $H_1$)  & $91.43\%$ & $88.56\%$ & $84.53\%$ & $81.89\%$ & $81.14\%$ & $77.75\%$ \\
$\Delta_2$ (from $H_2$)  & $86.33\%$ & $73.91\%$ & $73.62\%$ & $80.09\%$ & $77.12\%$ & $77.35\%$ \\
$\Delta_{\text{max}} = \max\pb{\Delta_0,\Delta_1,\Delta_2}$ & $95.72\%$ & $93.65\%$ & $94.44\%$ & $82.43\%$ & $78.80\%$ & $79.78\%$ \\
\hline
\end{tabu}
}
\end{table}
\vspace{-4mm}

\begin{figure}
  \centering
  \begin{subfigure}[b]{0.24\linewidth}
    \includegraphics[width=\linewidth]{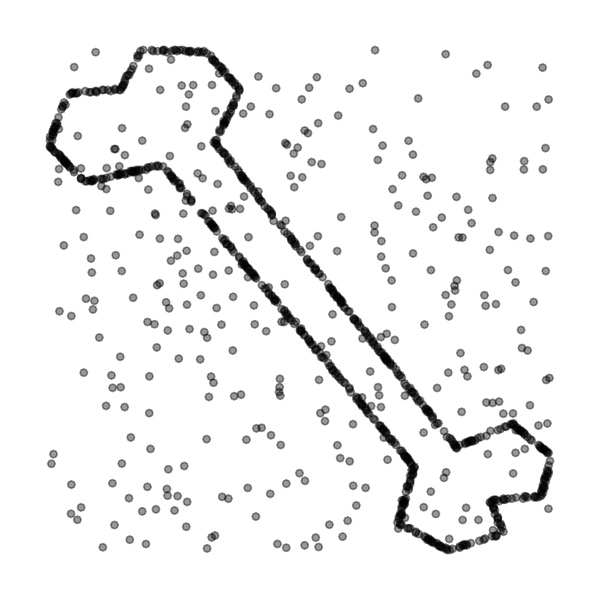}
    \caption{}
  \end{subfigure}
  \begin{subfigure}[b]{0.24\linewidth}
    \includegraphics[width=\linewidth]{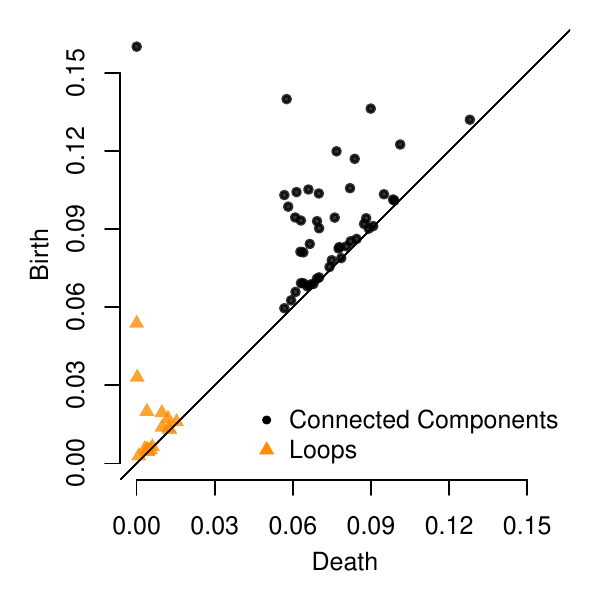}
    \caption{}
  \end{subfigure}
  \begin{subfigure}[b]{0.24\linewidth}
    \includegraphics[width=\linewidth]{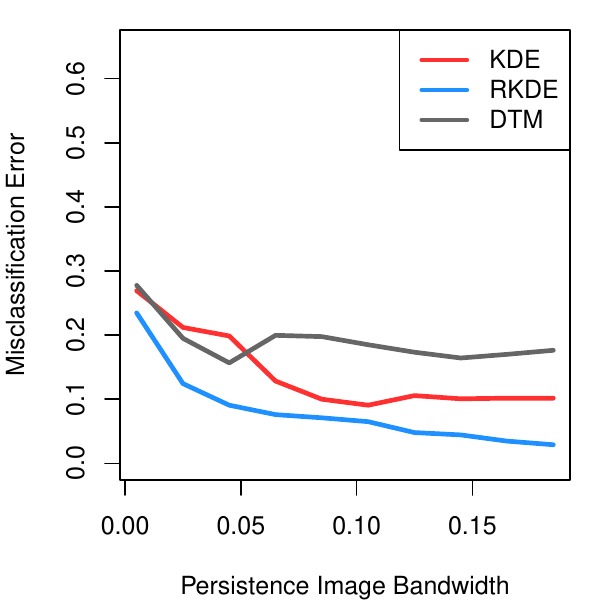}
    \caption{}
  \end{subfigure}
  \begin{subfigure}[b]{0.24\linewidth}
    \includegraphics[width=\linewidth]{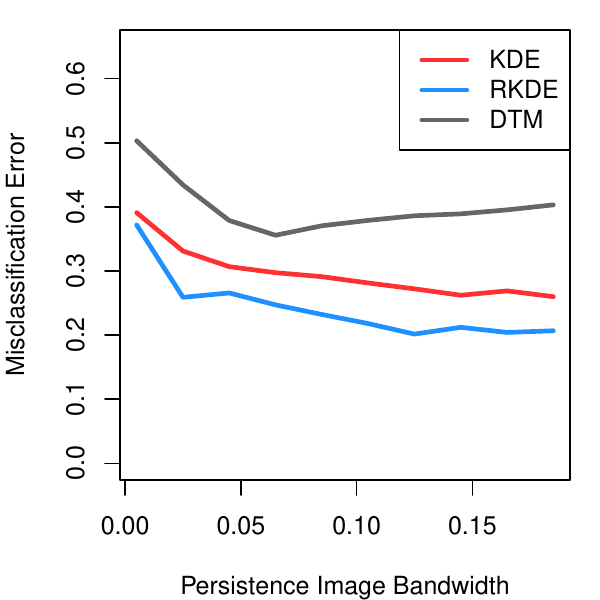}
    \caption{}
  \end{subfigure}
  \caption{(a) $\Xn$ is sampled from the image boundary of a \textit{bone}, and uniform noise is added. (b) The resulting persistence diagram from the robust KDE. The persistence diagram picks up the $1^{st}$--order features near the joints of the cartoon bone. The misclassification error for the KDE, robust KDE and DTM as the persistence image bandwidth increases, (c) for the three-class classification and, (d) for the five-class classification.}  \label{fig:mpeg7}
\end{figure}

\textbf{MPEG7.} In this experiment, we examine the performance of persistence diagrams in a classification task on \cite{latecki2000shape}. For simplicity, we only consider five classes: \textit{beetle}, \textit{bone}, \textit{spring}, \textit{deer} and \textit{horse}. We first extract the boundary of the images using a Laplace convolution, and sample $\Xn$ uniformly from the boundary of each image, adding uniform noise ($\pi=15\%$) in the enclosing region. Persistence diagrams $\dgm\pa{\barfn}$ and $\dgm\pa{\fn}$ from the KDE and robust KDE are constructed. In addition, owing to their ability to capture nuanced multi-scale features, we also construct $\dgm\pa{d_{n,m}}$ from the DTM filtration. The smoothing parameters $\s(k)$ and $m(k)$ are chosen as earlier for $k=5$. The persistence diagrams are normalized to have a max persistence $\max\{\abs{d-b} = 1:(b,d) \in \dgm(\phi)\}$, and then vectorized as persistence images, $\Img\pa{\barfn,h}$, $\Img\pa{\fn,h}$, and $\Img\pa{d_{n,m},h}$ for various bandwidths $h$. A linear SVM classifier is then trained on the resulting persistence images. In the first experiment we only consider the first three classes, and in the second experiment we consider all five classes. The results for the classification error, shown in Figure~\ref{fig:mpeg7}, demonstrate the superiority of the proposed method. We refer the reader to Appendix~\ref{additional-experiments} for 
additional experiments.

%


\section{Conclusion \& Discussion}
\label{discussion}

In this paper, we proposed a statistically consistent robust persistent diagram using RKHS-based robust KDE as the filter function. By generalizing the notion of influence function to the space of persistence diagrams, we mathematically and empirically demonstrated the robustness of the proposed method to that of persistence diagrams induced by other filter functions such as KDE. Through numerical experiments, we demonstrated the advantage of using robust persistence diagrams in machine learning applications. 
%
%
We would like to highlight that most of the theoretical results of this paper crucially hinge on the loss function being convex. As a future direction, we would like to generalize the current results to non-convex loss functions, and explore robust persistence diagrams induced other types of robust density estimators, which could potentially yield more robust persistence diagrams. {Another important direction we intend to explore is to enhance the computational efficiency of the proposed approach using coresets, as in \cite{brecheteau2018k}, and/or using weighted Rips filtrations, as in \cite{anai2019dtm}. We provide a brief discussion in Appendix~\ref{persistent-homology}.}

\section*{Broader Impact}
\label{impact}
Over the last decade, Topological Data Analysis has become an important tool for extracting geometric and topological information from data, and its applications have been far reaching. For example, it has been used successfully in the study the fragile X-syndrome, to discover traumatic brain injuries, and has also become an important tool in the study of protein structure. In astrophysics, it has aided the study of cosmic microwave background, and the discovery of cosmic voids and filamental structures in cosmological data. With a continual increase in its adoption in data analysis, it has become important to understand the limitations of using persistent homology in machine learning applications. As real-world data is often flustered with measurement errors and other forms of noise, in this work, we examine the sensitivity of persistence diagrams to such noise, and provide methods to mitigate the effect of this noise, so as to make reliable topological inference.

\begin{ack}
  {The authors would like to thank the anonymous reviewers for their helpful comments and constructive feedback. Siddharth Vishwanath and Bharath Sriperumbudur are supported in part by NSF DMS CAREER Award 1945396.
  Kenji Fukumizu is supported in part by JST CREST Grant Number JPMJCR15D3, Japan. Satoshi Kuriki is partially supported by JSPS KAKENHI Grant Number JP16H02792, Japan.}
\end{ack}




\renewcommand\refname{References}
\bibliographystyle{abbrvnat}

\bibliography{nips-refs}

\newpage
\setcounter{page}{1}
\appendix
\numberwithin{equation}{section}

\rule{\linewidth}{1pt}

\begin{center}
  \textbf{\large{Supplementary Material:}}%
  \\ \large{Robust Persistence Diagrams using Reproducing Kernels}
\end{center}

\rule{\linewidth}{2pt}


\section{Proofs for Section \ref{results}}
\label{proofs}

In what follows, given a metric space $\pa{\X,\varrho}$, $L_p(\X,\mu)$ is the Banach space of functions of $p^{th}$-power $\mu$-integrable functions with norm $\Vert \cdot\Vert_{p}$, where $\mu$ is a Borel measure defined on $\X$. For a fixed loss $\rho$, we will use the notation $\ell_g(\cdot) = \ell(\cdot,g) = \rho\pa{\norm{\Phi(\cdot) - g}_{\HH}}$ in order to emphasize the dependency of the loss on the choice of $g \in \G$. Borrowing some notation from empirical process theory, we define the empirical risk-functional in Eq.~\eqref{def:rkde} as
\eq{
\J_n(g) \defeq \pr_n\ell_g = \sum_{i=1}^{n}\rho\pa{\normh{\Phis(\Xv_i)-g}},\nn
}
and, similarly, the population risk functional $\J(g)$ is given by
\eq{
\J(g) \defeq \pr\ell_g = \int_{\R^d}{\rho\pa{\normh{\phigox}}} \ d\pr(\xv).\nn
}

\subsection{Proof of Theorem \ref{thm:influence}}
\label{proof:influence}

For $\e > 0$, define the risk functional associated with $\prx$ to be
\eq{
\Jx(g) = \prx\ell_g = (1-\e)\J(g) + \e\rho\pa{\normh{\phigox}},\nn
}
and let $\fox = \arginf_{g\in \G}\Jx(g)$ be its minimizer. From the stability result of Proposition \ref{prop:stability} we have that
\eq{
\Psi\pa{\fo;\xv} = \lim_{\epsilon \rightarrow 0}\f{1}{\epsilon} \Winf\pa{\dgm\pa{\fox},\dgm\pa{\fo}} \le \lim_{\epsilon \rightarrow 0}\f{1}{\epsilon} \norminf{\fox-\fo}.\nn
}
Using Propositions \ref{lemma:gamma} and \ref{lemma:equi}, we know that the sequence $\pb{\Jx}$ is equi-coercive, and \mbox{$\Jx$ ${\Gamma}$--converges to $\J$} as $\e \rightarrow 0$. From the fundamental theorem of $\Gamma$--convergence \cite[Theorem 7.8]{dal2012introduction} we have that ${\normh{\fox-\fo} \rightarrow 0}$, and, consequently, from Lemma \ref{lemma:supnorm}, $\norminf{\fox-\fo} \rightarrow 0$ as $\epsilon \rightarrow 0$. Thus,
\eq{
\lim_{\epsilon \rightarrow 0}\f{1}{\epsilon} \norminf{\fox-\fo} = \norminf{\lim_{\epsilon \rightarrow 0}\f{\fox-\fo}{\epsilon}}.
\label{eq:suplimits}
}

Let the limit in the right hand side of Eq.~\eqref{eq:suplimits} be denoted by $\dotfo$. Although $\dotfo$ does not admit a closed-form solution, from \cite[Theorem 8]{kim2012robust} we have that $\dotfo$ satisfies $V = a\dotfo + B$, where
\eq{
V &= \varphi\pa{\normh{\phifox}}\cdot\pa{\phifox},\nn\\
a &= {\int_{\R^d}{\varphi\pa{\normh{\phifo}}d\pr(\yv)}}, \ \ \ \ \text{and}\nn\\
B &= \int_{\R^d}\pa{\f{\varphi'\pa{\normh{\phifo}}}{\normh{\phifo}} \Innerh{\dotfo,\phifo} \cdot \pa{\phifo}}d\pr(\yv).\nn
}

For brevity, we adopt the notation $z(\yv) = \normh{\phifo}$ and $u(\cdot, \yv) = \f{\phifo}{\normh{\phifo}} \in \HH$. Then note that
$a \in \R$ and $B \in \HH$ are given by
\eq{
a &= \int_{\R^d}{\varphi\pa{z(\yv)}d\pr(\yv)}, \ \ \ \ \text{and} \nn\\
B &= \int_{\R^d}{z(\yv)\varphi'\pa{z(\yv)}\Innerh{\dotfo,u(\cdot, \yv)}  u(\cdot,\yv) \ d\pr(\yv)}.\nn
}
Using the reverse triangle inequality we have
\eq{
\normh{V} \ge a\normh{\dotfo} - \normh{B}.
\label{eq:normv}
}
We now look for an upper bound on $\normh{B}$. By noting that
\eq{
\Innerh{\dotfo,u(\cdot,\xv)}\Innerh{\dotfo,u(\cdot,\yv)}\Innerh{u(\cdot,\xv),u(\cdot,\yv)} \le \normh{\dotfo}^2,\nn
}
we have
\eq{
\normh{B}^2 = \Innerh{B,B} &\le \iint{z(\xv)\varphi'(z(\xv))z(\yv)\varphi'(z(\yv))\normh{\dotfo}^2 d \pr(\xv) \ d\pr(\yv)}\nn\\
&= \normh{\dotfo}^2 \pa{\int_{\R^d}{z(\yv)\varphi'(z(\yv)) \ d\pr(\yv)}}^2.\nn
}
Plugging this back into Eq.~\eqref{eq:normv} we get
\eq{
\normh{V} &\ge \normh{\dotfo} \int_{\R^d}{\varphi(z(\yv)) - z(\yv)\varphi'(z(\yv))} \ d\pr(\yv)\nn\\
&= \normh{\dotfo} \int_{\R^d}{\zeta(z(\yv)) \ d\pr(\yv)},
\label{eq:lower-normv}
}
where $\zeta(z) = \varphi(z) - z\varphi'(z)$. Similarly, by using the definition of $\varphi$, it follows that
\eq{
\normh{V} = \varphi\pa{\normh{\phifox}}\cdot\normh{\phifox} = \rho'\pa{\normh{\phifox}}.\nn
}
Combining this with Eq.~\eqref{eq:lower-normv} we get
\eq{
\normh{\dotfo} \le \f{\rho'\pa{\normh{\phifox}}}{\int_{\R^d}{\zeta\pa{\normh{\phifo}}}d\pr(\yv)}.\nn
}
By noting that $\norminf{\dotfo} \le \kinf^\half\normh{\dotfo}$ and $\Psi(\fo;\xv) \le \norminf{\dotfo}$, the result follows. \qed


\subsection{Proof for Theorem \ref{thm:consistency2}}
\label{proof:consistency2}

Using the triangle inequality we can break our problem down as follows
\eq{
\norminf{\fn-f} \le \underbrace{\norminf{\barfo-f}}_{\circled{a}} + \underbrace{\norminf{\fo-\barfo}}_{\circled{b}},\nn
}
where, $\barfo = \int_{\R^d}{\K(\cdot,\xv)d\pr(\xv)}$ is the population level KDE. For term $\circled{a}$, for $\pr \in \mathcal{M}(\R^d)$, it is well known \citep{chen2017tutorial} that the approximation error for the KDE vanishes, i.e.,
\eq{
\norminf{\barfo-f} \rightarrow 0,\nn
}
as $\sigma \rightarrow 0$. So, it remains to verify that $\circled{b}$ vanishes, i.e., $\norminf{\fo-\barfo} \rightarrow 0$. With this in mind, consider the map $\Ts : \G \rightarrow \G$ given by
\eq{
\Ts(g) = \int\limits_{\R^d}{\f{\varphi\pa{\normh{\Phis\pa{\xv}-g}}}{\int\limits_{\R^d}\varphi\pa{\normh{\Phis\pa{\xv}-g}}d\pr(\xv)}} \Phis(\xv)d\pr(\xv).\nn
}
Our approach to verifying that $\circled{b}$ vanishes is similar to \citet[Lemma 9]{vandermeulen2013consistency}, where we show that the map $\Ts$ is a contraction map when restricted to the subspace
\eq{
\Qs \defeq B_{\HH}\pa{\zerov,\delta\nus} \cap \D.\nn
}
A key difference is that we work with $\norminf{\cdot}$--norm, requiring us to obtain a sharper bound for the Lipschitz constant associated with the contraction.

For brevity, we adopt the notation $m(\xv,g) = \varphi\pa{\normh{\Phis\pa{\xv}-g}}$. The authors in \cite{kim2012robust} show that $\fo$ is a fixed point of the map $\Ts$, i.e., $\Ts(\fo) = \fo$, and that $\barfo$ is the image of $\zerov$ under $\Ts$, i.e., $\Ts(\zerov) = \barfo$. Additionally, from Lemma \ref{lemma:vanishing}, we know that $\normh{\fs} \le \delta \nus$, for some $0 < \delta < 1$. Thus, we can rewrite $\fo-\barfo = \Ts(\fo) - \Ts(\zerov)$.

Let $g,h \in \Qs$. Then we have that
\eq{
\Ts(g) - \Ts(h) &= \int\limits_{\R^d}{\f{m(\xv,g)}{\int\limits_{\R^d}m(\yv,g)d\pr(\yv)}} \Phis(\xv)d\pr(\xv) - \int\limits_{\R^d}{\f{m(\uv,h)}{\int\limits_{\R^d}m(\vv,h)d\pr(\vv)}} \Phis(\uv)d\pr(\uv)\nn\\
&= \f{1}{\alpha\beta}\cdot \pa{%
\beta \int\limits_{\R^d}{m(\xv,g)\Phis(\xv)d\pr(\xv)} -  %
\alpha \int\limits_{\R^d}{m(\uv,h)\Phis(\uv)d\pr(\xv)}%
}\nn\\
&= \f{1}{\alpha\beta}\cdot \xi,
\label{eq:ts-def}
}
where $\alpha \defeq \int_{\R^d}m(\yv,g)d\pr(\yv) \in \R$, $\beta \defeq \int_{\R^d}m(\vv,h)d\pr(\vv) \in \R$ and the numerator $\xi \in \HH$.

By Tonelli's theorem
\eq{
\xi &= \beta \int\limits_{\R^d}{m(\xv,g)\Phis(\xv)d\pr(\xv)} -  %
\alpha \int\limits_{\R^d}{m(\uv,h)\Phis(\uv)d\pr(\xv)}\nn\\
&= \int\limits_{\R^d}{m(\xv,g)\Phis(\xv)\pa{{\int\limits_{\R^d}m(\vv,h)d\pr(\vv)}}d\pr(\xv)}\nn\\   %
&\qquad\qquad-\int\limits_{\R^d}{m(\uv,h)\Phis(\uv)\pa{{\int\limits_{\R^d}m(\yv,g)d\pr(\yv)}}d\pr(\xv)}\nn\\
&= \iint\limits_{\R^d\times\R^d}m(\xv,g)m(\vv,h)\Phis(\xv) d\pr(\vv)d\pr(\xv) - %
\iint\limits_{\R^d\times\R^d}m(\uv,h)m(\yv,g)\Phis(\uv) d\pr(\yv)d\pr(\uv)\nn\\
&= \iint\limits_{\R^d\times\R^d}{\Phis(\xv) \pc{m(\xv,g)m(\yv,h) - m(\xv,h)m(\yv,g)} d\pr(\xv)d\pr(\yv)}.
\label{eq:xi}
}

Then by adding and subtracting $m(\xv,h)m(\yv,h)$ to the term inside, we get
\eq{
m(\xv,g)m(\yv,h) - m(\xv,h)m(\yv,g) &= m(\yv,h)\pb{m(\xv,g)-m(\xv,h)}\nn\\
&\qquad\qquad+ m(\xv,h)\pb{m(\yv,h)-m(\yv,g)}.\nn
}

Plugging this back into Eq.~\eqref{eq:xi}, we get $\xi = \xi_1 + \xi_2$ where
\eq{
\xi_1 &= \iint\limits_{\R^d\times\R^d}\Phis(\xv)\pb{m(\xv,g)-m(\xv,h)}m(\yv,h)d\pr(\yv)d\pr(\xv)\nn\\
&= \int\limits_{\R^d}{m(\yv,h)d\pr(\yv)}\int\limits_{\R^d}{\Phis(\xv)\pb{m(\xv,g)-m(\xv,h)} d\pr(\xv)}\nn\\
&= \beta \int\limits_{\R^d}{\K(\cdot,\xv)\pb{m(\xv,g)-m(\xv,h)} d\pr(\xv)}\nn\\
&\stackrel{(i)}{=} \beta \cdot \pc{\psis*\pa{\pa{m(\cdot,g)-m(\cdot,h)}f(\cdot)}},\nn
}
where (i) follows from the fact that the kernel $\K(\xv,\yv) = \psis(\xv-\yv)\defeq \sigma^{-d}\psi(\Vert \xv-\yv\Vert_2/\sigma)$ is translation invariant and $f$ is the density associated with $\pr$.

Similarly,
\eq{
\xi_2 &= \iint\limits_{\R^d\times\R^d}\Phis(\xv)m(\xv,h)\pb{m(\yv,h)-m(\yv,g)}d\pr(\xv)d\pr(\yv)\nn\\
&= \int\limits_{\R^d}{\pc{m(\yv,h)-m(\yv,g)} d\pr(\yv)}\int\limits_{\R^d}{\Phis(\xv)m(\xv,h)d\pr(\xv)}\nn\\
&\stackrel{}{\le} \norminf{m(\cdot,h)-m(\cdot,g)} \cdot \pc{\psis*\pa{m(\cdot,h)f(\cdot)}}.\nn
}
The upper bound for $\norminf{\xi_1}$ is as follows
\eq{
\norminf{\xi_1} &= \beta \norminf{\psis*\pa{\pa{m(\cdot,g)-m(\cdot,h)}f(\cdot)}}\nn\\
&\stackrel{(i)}{\le} \beta \norm{\psis}_1 \norminf{\pa{m(\cdot,g)-m(\cdot,h)}f(\cdot)}\nn\\
&\stackrel{(ii)}{\le} \beta \norminf{m(\cdot,g)-m(\cdot,h)}\norminf{f},
\label{eq:xi1}
}
where (i) follows from Young's inequality \citep[Theorem 20.18]{hewitt1979abstract} and (ii) follows from the fact that $\norminf{fg} \le \norminf{f}\norminf{g}$. Similarly, for $\xi_2$ we have
\eq{
\norminf{\xi_2} &\le \norminf{m(\cdot,h)-m(\cdot,g)} \norminf{\psis*\pa{m(\cdot,h)f(\cdot)}}\nn\\
&\stackrel{(i)}{\le} \norminf{m(\cdot,h)-m(\cdot,g)} \norm{\psis}_1 \norminf{m(\cdot,h)f(\cdot)}\nn\\
&\stackrel{(ii)}{\le} \norminf{m(\cdot,h)-m(\cdot,g)} \norminf{m(\cdot,h)} \norminf{f}.
\label{eq:xi2}
}

From the proof of \cite[Lemma 9, Page 20--22]{vandermeulen2013consistency}, for $g,h \in \Qs$ for fixed constants $c_1,c_2 >0$ we have the following two bounds:
\eq{
\alpha,\beta \ge \f{1}{c_1\nus},
\label{eq:xi-temp1}
}
and
\eq{
\norminf{m(\cdot,h)-m(\cdot,g)} \le \normh{g-h} c_2\nus^{-2},
\label{eq:xi-temp2}
}
where the last inequality follows from the Lipschitz property of $\varphi$ and fact that $\rho$ is strictly convex. Additionally, for $c_3 = \norminf{\rho'} < \infty$ we have
\eq{
m(\xv,g) &= \varphi\pa{\normh{\Phis\pa{\xv}-g}}\nn\\
&= \f{\rho'\pa{\normh{\Phis\pa{\xv}-g}}}{{\normh{\Phis\pa{\xv}-g}}}\nn\\
&\le \f{c_3}{{\normh{\Phis\pa{\xv}-g}}}\nn\\
&\stackrel{(iii)}{\le} \f{c_3}{\Abs{\normh{\Phis(\xv)} - \normh{g}}}\nn\\
&= \f{c_3}{(1-\delta)\nus},
\label{eq:xi-temp3}
}
where (iii) follows from reverse triangle inequality. Plugging the bounds in equations \eqref{eq:xi-temp1}, \eqref{eq:xi-temp2} and \eqref{eq:xi-temp3} back into equations \eqref{eq:xi1} and \eqref{eq:xi2} we get,
\eq{
\norminf{\xi_1}+\norminf{\xi_2} \le \norminf{f} \pa{{\beta c_2}\nus^{-2}\normh{g-h} + \f{c_2c_3}{(1-\delta)}\nus^{-3}\normh{g-h}}.\nn
}
Using this upper bound in Eq.~\eqref{eq:ts-def} we get
\eq{
\norminf{\Ts(g) - \Ts(h)} &= \norminf{\f{\xi}{\alpha\beta}}\nn\\
&\le \f{\norminf{\xi_1}+\norminf{\xi_2}}{\alpha\beta}\nn\\
&\stackrel{(iv)}{\le} \norminf{f}\pa{\f{c_1c_2}{c_1}\nus^{-1}\normh{g-h} + \f{c_2c_3}{c_1^2(1-\delta)}\nus^{-1}\normh{g-h}}\nn\\
&\stackrel{(v)}{=} C\nus^{-1}\normh{g-h}\nn\\
&\stackrel{(vi)}{\le} C\nus\inv\norminf{g-h}^{\half},\nn
}
where in (iv) we use Eq.~\eqref{eq:xi-temp1}, in (v) we use the fact that whenever $\pr \in \M(\R^d)$, we have $\norminf{f} < \infty$ and $C > 0$ is a constant depending only on $c_1,c_2,c_3$ and $\norminf{f}$. Additionally, (vi) holds through an application of Lemma \ref{lemma:supnorm} to $g-h \in \Qs \subset \D$. This confirms that $\Ts$ is a contraction mapping. We use this to show that $\circled{b}$ vanishes as $\sigma \rightarrow 0$. Since $\fo,\zerov \in \Qs$ and $\fo-\zerov \in \D$, we have that
\eq{
\norminf{\fo-\barfo} &= \norminf{\Ts(\fo)-\Ts(\zerov)}\nn\\
&\le C\nus\inv \norminf{\fo-\zerov}^\half\nn\\
&= C\nus\inv \norminf{\fo}^\half.\nn
}
Using the triangle inequality $\norminf{\fo}^\half \le \norminf{\fo-\barfo}^\half + \norminf{\barfo}^\half$ we get
\eq{
\norminf{\fo-\barfo} &\le C\nus\inv \pa{\norminf{\fo-\barfo}^\half + \norminf{\barfo}^\half}\nn\\
&=C\nus\inv \pa{\norminf{\Ts\pa{\fo}-\Ts(\zerov)}^\half + \norminf{\barfo}^\half}\nn\\
&\le C\nus\inv\pa{\pa{C\nus\inv\norminf{\fo-\zerov}^\half}^\half + \norminf{\barfo}^\half}\nn\\
& = C^{\f{3}{2}}\nus^{-\f{3}{2}}\norminf{\fo}^{\f{1}{4}} + C\nus\inv\norminf{\barfo}^\half,
\label{eq:contraction-main}
}
by using the contraction mapping twice. Observe that
\eq{
\norminf{\fo} &\le \nus \normh{\fo} \le \delta\nus^2,\nn
}
where the first inequality follows from Lemma \ref{lemma:supnorm} and the second inequality follows from the fact that $\normh{\fo} \le \delta\nus$ since $\fo \in \Qs$. Furthermore, $\norminf{\barfo} = \norminf{\psis*f} \le \norminf{\psis} \norm{f}_1 \le \nus$ from Young's inequality. By noting that $\nus = \psis(0) = \s^{-d}\psi(0)$, collecting these bounds back into Eq.~\eqref{eq:contraction-main} we get
\eq{
\norminf{\fo-\barfo} \le C^{\f{3}{2}}\delta^{\f{1}{4}}\nus\inv + C\nus^{-\half}\sqrt{\psi(0)}.\nn
}

yielding that $\norminf{\fo-\barfo} \rightarrow 0$ as $\sigma \rightarrow 0$, thereby verifying that $\circled{b}$ vanishes as $\sigma \rightarrow 0$.
\qed


\subsection{Proof of Theorem \ref{thm:consistency}}
\label{proof:consistency}

The proof proceeds in two steps: We first establish the uniform consistency for the robust KDE and then use the bottleneck stability to show consistency of the robust persistence diagrams in $\Winf$. From the stability theorem for persistence diagrams \citep{cohen2007stability,chazal2016structure}, we have that $\Winf\pa{\dgm\pa{\fn},\dgm\pa{\fo}} \le \norm{\fn-\fo}_\infty$. Thus, it suffices to show that $\norm{\fn-\fo}_\infty \stackrel{p}{\rightarrow} 0$ as $n \rightarrow \infty$. In order to prove the latter, we adapt the argmax consistency theorem \citep[Theorem 5.7]{van2000asymptotic} for minimizers of a risk function.

\pagebreak

\begin{lemma}[Theorem 5.7, \cite{van2000asymptotic}]
  Given a metric space $\pa{\G,d}$, let $\mathcal{J}_n$ be random functions and $\mathcal{J}$ be a fixed function of $g \in \G$ such that for every $\epsilon > 0$,
  \begin{enumerate}[label=(\arabic*)]
    \item $\inf\limits_{g: d(g,g_0) \ge \epsilon}\mathcal{J}(g) >  \mathcal{J}(g_0)$, and
    \item $\sup\limits_{g\in\G}\abs{\mathcal{J}_n(g) - \mathcal{J}(g)} \stackrel{p}{\rightarrow} 0$.
  \end{enumerate}
  Then any sequence $g_n$ satisfying $\mathcal{J}_n(g_n) < \mathcal{J}_n(g_0) + O_p(1)$ satisfies $d(g_n,g_0) \stackrel{p}{\rightarrow} 0$.
  \label{lemma:argmaxconsistency}
\end{lemma}

For $\G = \HH \cap \D$, and $d(\fn,\fo) = \norm{\fn-\fo}_\infty$, in order to establish uniform consistency of the robust KDE, as per Lemma \ref{lemma:argmaxconsistency}, we need to verify that conditions (1) and (2) are satisfied.

Condition (1) follows from the strict convexity of $\J(g)$ in Proposition \ref{lemma:lipschitz}. Specifically, \cite{kim2012robust} establish that assumptions $\pa{\mathcal{A}1} - \pa{\mathcal{A}3}$ guarantee the existence and uniqueness of $\fo = \arginf_{g\in\G}\J(g)$. Then, for any $g \in \G$ such that $\norm{g-\fo}_\HH > \delta$, we have that $\J(g) > \J(\fo)$.

We now turn to verifying condition (2). Observe that $\sup_{g\in\G}\abs{\mathcal{J}_n(g) - \mathcal{J}(g)}$ can be rewritten as the supremum of an empirical process, i.e.,
\eq{
\sup\limits_{g\in\G}\abs{\J_n(g) - \J(g)} = \sup\limits_{\ell_g \in \widetilde{\F}}\abs{\pr_n \ell_g - \pr \ell_g} \defeq \norm{\pr_n - \pr}_{\widetilde{\F}},\nn
}
where $\widetilde{\F} = \pb{\ell_g : g \in \G}$, and $\ell_g(\xv) = \rho\pa{\norm{\Phis(\xv)-g}_\HH}$. Verifying condition (2) reduces to showing that $\widetilde{\F}$ is a Glivenko-Cantelli class. Define $\eta(\cdot) = \normh{\Phis(\cdot) - g}^2$ and let $\F = \pb{\eta_g : g \in \G}$. For the continuous map $\xi: [0,\infty) \rightarrow [0,\infty)$ given by $\xi(t) = \rho(\sqrt{t})$, we have that
\eq{
\xi\circ \F = \pb{\xi(f) : f \in \F} = \pb{\xi\circ\eta_g(\cdot) : g \in \G} = \pb{\rho(\norm{\Phis(\cdot) - g}_\HH) : g \in \G}  = \widetilde{\F}.\nn
}
By the preservation theorem for Glivenko-Cantelli classes \citep[Theorem 3]{van2000preservation}, it holds that if $\F$ is a Glivenko-Cantelli class, then $\widetilde{\F}$ is also a Glivenko-Cantelli class. So verifying condition (2) reduces to verifying that $\F$ is a Glivenko-Cantelli class. To this end, we first show that $F(\xv_{1:n})=F(\xv_1,\xv_2,\dots,\xv_n) = \sup_{g \in \G}\abs{\pr_n\eta_g - \pr\eta_g} = \norm{\pr_n - \pr}_{\F}$ satisfies the self-bounded property for McDiarmid's inequality, i.e.,
\eq{
\sup\limits_{\xv_i\neq\xv'_i}\abs{F(\xv_{1:n}) - F(\xv'_{1:n})} &\le \f{1}{n} \sup_{\xv_i,\xv'_i}\sup_{g \in \G}\pa{ \norm{\Phis(\xv_i)}^2_\HH + \norm{\Phis(\xv'_i)}^2_\HH + 2\abs{g(\xv_i)} + 2\abs{g(\xv'_i)}}.\nn
}
Observe that $\norm{\Phis(\xv)}^2_\HH  = \K(\xv,\xv) \le \norm{\K}_\infty$ and $\abs{g(\xv)} \le \norm{g}_\infty < \norm{\K}_\infty$ by Lemma \ref{lemma:supnorm}. Thus, we have that
\eq{
\sup\limits_{\xv_i\neq\xv'_i}\abs{F(\xv_{1:n}) - F(\xv'_{1:n})} \le \f{6\kinf}{n}.\nn
}
From \cite[Theorem 9]{bartlett2002rademacher}, we have that with probability greater than $1 - e^{-\delta}$,
\eq{
\norm{\pr_n - \pr}_{\F} \le 2 \mathfrak{R}_n(\F) + \sqrt{\f{3\delta\kinf}{n}},
\label{eq:rademacher-consistency}
}

where $\mathfrak{R}_n(\F)$ is the Rademacher complexity of $\F$ given by,
\eq{
\mathfrak{R}_n(\F) &= \E_{\epsilon}\pa{\sup_{g\in\G} \abs{\f{1}{n}\sum_{i=1}^n {\epsilon_i \norm{\Phis(\xv_i)-g}^2_\HH} }}\nn\\
&\le \E_{\epsilon}\pa{\sup_{g\in\G} \pb{\abs{\f{1}{n}\sum_{i=1}^n {\epsilon_i \norm{\Phis(\xv_i)}^2_\HH}} + \abs{\f{1}{n}\sum_{i=1}^n {\epsilon_i \norm{g}^2_\HH}} + 2\abs{\f{1}{n}\sum_{i=1}^n {\epsilon_i g(\xv_i)}}} }\nn\\
&= \circled{1} + \circled{2} + \circled{3}.\nn
}
Note that $\E_{\epsilon}\pa{f(\epsilon_{1:n},\xv_{1:n})} \defeq \E\pa{f(\epsilon_{1:n},\xv_{1:n}) | \xv_{1:n}}$ is the conditional expectation of the Rademacher random variables $\epsilon_1,\epsilon_2,\dots,\epsilon_n$, keeping $\xv_1,\xv_2,\dots,\xv_n$ fixed. First, we have that,
\eq{
\circled{1} = \E_{\epsilon} \pa{\sup_{g\in\G} \abs{\f{1}{n}\sum_{i=1}^n {\epsilon_i \norm{\Phis(\xv_i)}^2_\HH} }} &\stackrel{(i)}{=} \E_{\epsilon} {\abs{\f{1}{n}\sum_{i=1}^n {\epsilon_i \K(\xv_i,\xv_i)} }}\nn\\
&\stackrel{(ii)}{\le} \sqrt{\E_{\epsilon}{\abs{\f{1}{n}\sum_{i=1}^n {\epsilon_i \K(\xv_i,\xv_i)} }^2}}\nn\\
&\le \sqrt{\E_{\epsilon} \pa{\f{1}{n^2}\sum_{i,j}\epsilon_i\epsilon_j\K(\xv_i,\xv_i)\K(\xv_j,\xv_j)}}\nn\\
&\stackrel{(iii)}{=} \f{1}{\sqrt{n}}\kinf,\nn
}
where (i) follows from the absence of $g$ inside the expectation, (ii) follows from Jensen's inequality and (iii) follows from the fact that $\epsilon_i \independent \epsilon_j$ for $i\neq j$. For the second term, we have
\eq{
\circled{2} = \E_{\epsilon} \pa{\sup_{g\in\G} \abs{\f{1}{n}\sum_{i=1}^n {\epsilon_i \norm{g}^2_\HH} }} &= \E_{\epsilon} \pa{\sup_{g\in\G} \norm{g}^2_\HH \abs{\f{1}{n}\sum_{i=1}^n {\epsilon_i} }}\nn\\
&\le \sup_{g\in\G} \norm{g}^2_\HH \sqrt{\E_{\epsilon}{\abs{\f{1}{n}\sum_{i=1}^n {\epsilon_i}}^2}},\nn\\
&\stackrel{(iv)}{\le} \f{1}{\sqrt{n}}\kinf,\nn
}
where (iv) follows from the fact that $\norm{g}^2_\HH \le \kinf$. Lastly, we have
\eq{
\circled{3} = 2 \E_{\epsilon} \pa{\sup_{g\in\G} \abs{\f{1}{n}\sum_{i=1}^n {\epsilon_i g(\xv_i)} }} &\stackrel{(v)}{=} 2 \E_{\epsilon} \pa{\sup_{g\in\G} \abs{\Innerh{g,\f{1}{n}\sum_{i=1}^n {\epsilon_i \K(\cdot,\xv_i) }} }}\nn\\
&\stackrel{(vi)}{\le} 2 \E_{\epsilon} \pa{\sup_{g\in\G} {\normh{g}\normh{\f{1}{n}\sum_{i=1}^n {\epsilon_i \K(\cdot,\xv_i) }} }}\nn\\
&= 2 \sup_{g\in\G} \norm{g}_\HH \E_{\epsilon} \pa{\sqrt{\f{1}{n^2}\sum_{i,j}{\epsilon_i\epsilon_j\K(\xv_i,\xv_j)}}}\nn\\
&\stackrel{(vii)}{\le} 2 \f{\kinf^{\half}}{n} \sqrt{\E_{\epsilon} \pa{\sum_{i,j}{\epsilon_i\epsilon_j\K(\xv_i,\xv_j)}}}\nn\\
&\stackrel{(viii)}{\le} \f{2}{\sqrt{n}} \kinf,\nn
}
where (v) follows from the reproducing property, (vi) is obtained from Cauchy-Schwarz inequality, (vii) follows from Jensen's inequality, and (viii) follows from the fact that $\epsilon_i \independent \epsilon_j$ for $i\neq j$. Collecting these three inequalities, we have
\eq{
\mathfrak{R}_n(\F)  = \circled{1} + \circled{2} + \circled{3} \le \f{4}{ \sqrt{n}} \kinf.\nn
}

Plugging this into Eq.~\eqref{eq:rademacher-consistency}, we have with probability greater than $1 - e^{-\delta}$,
\eq{
\norm{\pr_n - \pr}_{\F} \le \f{8\kinf}{\sqrt{n}} + \sqrt{\f{3\delta\kinf}{n}},\nn
}
which implies that $\norm{\pr_n - \pr}_{\F} \rightarrow 0$ as $n \rightarrow \infty$, implying that $\F$ is a Glivenko-Cantelli class. The result, therefore, follows from Lemma \ref{lemma:argmaxconsistency}. \qed


\subsection{Proof of Theorem \ref{thm:concentration1}}
\label{proof:concentration}

For $g\in\G$ define the random fluctuation w.r.t. $\fo$ as
\eq{
\Delta\pa{\Xv,g} = \pa{\ell_g(\Xv) - \ell_{\fo}(\Xv)} - \pa{\J(g) - \J(\fo)}.\nn
}
The fluctuation process is an empirical process defined as
\eq{
\Delta_n(g) = \pr_n\Delta(\Xv,g) &= \pa{\J_n(g) - \J_n(\fo)} - \pa{\J(g) - \J(\fo)},\nn\\
&= \pr_n\pa{\ell_g-\ell_{\fo}} - \pr\pa{\ell_g-\ell_{\fo}}.\nn
}

We first show that the behaviour of $\normh{\fn -\fo}$ is controlled by the tail behaviour of the supremum of the fluctuation process. To this end, for $\delta > 0$, let
\eq{
\G_\delta = \pb{g \in \G : \normh{g-\fo} \le \delta} = B_{\HH}\pa{\fo,\delta} \cap \D.\nn
}
Suppose $\fn$ is such that $\normh{\fn - \fo} > \delta$, then, for sufficiently small $\lambda\in (0,1)$ such that $g=\lambda\fn + (1-\lambda)\fo \in \gdelta$, we have that
\eq{
\J_n(g) - \J_n(\fo) &\stackrel{(i)}{<} \lambda\J_n(\fn) + (1-\lambda)\J_n(\fo) - \J_n(\fo)\nn\\
&= \lambda \cdot \pa{\J_n(\fn) - \J_n(\fo)} \stackrel{(ii)}{\le} 0,
\label{eq:temp2}
}
where (i) follows from the strict convexity of $\J_n$ (Proposition \ref{lemma:lipschitz}), and (ii) follows from the fact that $\fn = \arginf_{g\in\G}\J_n(g)$. From Proposition \ref{lemma:lipschitz}, we also know that $\J$ is strongly convex such that
\eq{
\J(g) - \J(\fo) \ge \f{\mu}{2}\normh{g-\fo}^2.
\label{eq:temp3}
}
Combining equations \eqref{eq:temp2} and \eqref{eq:temp3} we have
\eq{
\f{\mu}{2}\normh{g-\fo}^2 &\le \J(g) - \J(\fo),\nn\\
&{=}  -\Bigl\{\pa{\J_n(g) - \J_n(\fo)} - \pa{\J(g) - \J(\fo)}\Bigr\} + \pa{\J_n(g) - \J_n(\fo)}\nn\\
&\le - \Delta_n(g) \le \supgdelta.\nn
}
By taking the supremum of the left hand side in the above inequality over all $g \in \gdelta$ we have
\eq{
\supgdelta \ge \f{\mu}{2}\delta^2
\label{eq:gdel}
}
This implies that whenever $\normh{\fn - \fo} > \delta$ holds, then the condition in Eq.~\eqref{eq:gdel} holds. Therefore,
\eq{
\pr^{\otimes n}\pb{\Xv_{1:n} : \normh{\fn - \fo} > \delta} \le \pr^{\otimes n}\pb{ \Xv_{1:n} : \supgdelta \ge \f{\mu}{2}\delta^2}.
\label{eq:pr-gdel}
}

We now study the behaviour of the r.h.s.~in Eq.~\eqref{eq:pr-gdel} using tools from empirical process theory. First, we show that $F(\xv_{1:n}) = F(\xv_1,\xv_2,\dots,\xv_n) = \sup\limits_{g\in\gdelta}\abs{\Delta_n(g)}$ satisfies the self-bounding property.
\eq{
\sup\limits_{\xv_i\neq\xv'_i}\abs{F(\xv_{1:n}) - F(\xv'_{1:n})} &= \sup\limits_{\xv_i\neq\xv'_i}\Abs{\sup\limits_{g\in\gdelta}\abs{\Delta_n(g)} - \sup\limits_{g\in\gdelta}\abs{\Delta_n(g)}},\nn\\
&\le \sup\limits_{\xv_i\neq\xv'_i}\sup\limits_{g \in \gdelta}\Abs{\Delta_n(g) - \Delta'_n(g)},\nn\\
&= \f{1}{n}\sup\limits_{\xv_i\neq\xv'_i}\sup\limits_{g \in \gdelta}\Abs{\pa{\ell_g(\xv_i) - \ell_{\fo}(\xv_i)} - \pa{\ell_g(\xv'_i) - \ell_{\fo}(\xv'_i)}},\nn\\
&\le \f{1}{n}\sup\limits_{\xv_i\neq\xv'_i}\sup\limits_{g \in \gdelta}\Abs{\pa{\ell_g(\xv_i) - \ell_{\fo}(\xv_i)}} + \Abs{\pa{\ell_g(\xv'_i) - \ell_{\fo}(\xv'_i)}},\nn\\
&\stackrel{(i)}{\le}\f{1}{n}\sup\limits_{g \in \gdelta}{2M\normh{g-\fo}} = \f{2M\delta}{n},\nn
}
where (i) follows from Proposition \ref{lemma:lipschitz} that $\ell_g$ is $M$-Lipschitz w.r.t. $\normh{\cdot}$. Therefore, from McDiarmid's inequality \citep[Theorem 2.9.1]{vershynin_2018} we have
\eq{
\pr^{\otimes n}\pb{\Xv_{1:n} : \supgdelta > \E\supgdelta + \epsilon} \le \exp\pa{-\f{n\epsilon^2}{2M^2\delta^2}}.
\label{eq:mcdiarmid1}
}
Next, we find an upper bound for the expected supremum of the fluctuation process. In order to do so, we first show that $\Delta_n(g)$ has sub-Gaussian increments. For fixed $g,h \in \G$ we have that $\E\pa{\Delta(\Xv,g)-\Delta(\Xv,h)} = 0$ and
\eq{
\Abs{\Delta(\Xv,g)-\Delta(\Xv,h)} \le \Abs{\ell_g(\Xv) - \ell_{h}(\Xv)} - \Abs{\J(g) - \J(h)} \le 2M \normh{g-h}.\nn
}
Since $\Abs{\Delta(\Xv,g)-\Delta(\Xv,h)}$ is bounded, it is, therefore, sub-Gaussian and from \citet[Example 2.5.8]{vershynin_2018}, we have that the sub-Gaussian norm $\norm{\Delta(\Xv,g)-\Delta(\Xv,h)}_{\psi_2} \le 2cM \normh{g-h}$ for $c>1/\sqrt{\log{2}}$. Consequently, the fluctuation process has sub-Gaussian increments with respect to the metric $\normh{g-h}$, i.e.,
\eq{
\norm{\Delta_n(g) - \Delta_n(h)}_{\psi_2} \le \f{1}{n} \sqrt{\sum_{i=1}^n{\norm{\Delta(\Xv_i,g) - \Delta(\Xv_i,h)}^2_{\psi_2}}}\le \f{M}{\sqrt{n}}\normh{g-h}.\nn
}

From the generalized entropy integral \cite[Lemma A.3]{srebro2010smoothness}, for a fixed constant $\gamma > 12/\sqrt{\log{2}}$ we have
\eq{
\E\supgdelta \le \inf\limits_{\alpha > 0}\pb{2\alpha + \f{\gamma M}{\sqrt{n}} \int_{\alpha}^{\delta}{\sqrt{\log{\mathcal{N}\pa{\gdelta,\normh{\cdot},\e}}}d\e }},
\label{eq:entropy1}
}
where $\mathcal{N}\pa{\gdelta,d,\e}$ is the $\e$-covering number of the class $\gdelta$ with respect to metric $d$.

We now turn our attention to finding an upper bound for $\mathcal{N}\pa{\gdelta,d,\e}$. Note that if $\B_\HH$ is a unit ball in the RKHS, then
\eq{
\log\mathcal{N}\pa{\gdelta,\normh{\cdot},\e} &= \log\mathcal{N}\pa{\B_\HH\cap\D,\normh{\cdot},\f{\e}{\delta}}\nn\\
&\stackrel{(i)}{\le} \log\mathcal{N}\pa{\B_\HH\cap\D,\norm{\cdot}_\infty,\pa{\f{\e}{\delta}}^2}\nn\\
&{\le} \log\mathcal{N}\pa{\B_\HH,\norm{\cdot}_\infty,\pa{\f{\e}{\delta}}^2},\nn
}
where (i) follows from Lemma \ref{lemma:supnorm} that $\normh{g-h}^2 \le \norminf{g-h}$. When the entropy numbers $e_n\pa{\id : \HH \rightarrow L_{\infty}(\X)}$ satisfy the assumption, from \cite[Lemma 6.21]{steinwart2008support} we have
\eq{
\log\mathcal{N}\pa{\B_\HH,\norm{\cdot}_\infty,\pa{\f{\e}{\delta}}^2} \le \pa{\f{a_\s\delta^2}{\e^2}}^{2p}.\nn
}

Plugging this into Eq.~\eqref{eq:entropy1}, we have that
\eq{
\E\supgdelta \le \inf\limits_{\alpha > 0}\pb{2\alpha + \f{\gamma Ma_\s\delta^{2p}}{\sqrt{n}} \int_{\alpha}^{\delta}{ \e^{-2p} d\e }} = \inf\limits_{\alpha > 0}T(\alpha),\nn
}
where $T(\alpha)$ is given by
\eq{
T(\alpha) =
\begin{cases}
  2\alpha + \gamma M\delta\sqrt{\f{a_\s}{n}}\log\pa{\f{\delta}{\alpha}} \ \ &\text{ if } p=\half,\nn\\ \\
  2\alpha + \f{\gamma M}{(1-2p)\sqrt{n}}\pa{\delta - \delta^{2p}\alpha^{1-2p}} \ \ &\text{ if } 0<p\neq\half<1.\nn
\end{cases}
}
At the value $\alpha_0$ where  $T(\alpha_0) = \inf_{\alpha>0}T(\alpha)$, we have
\eq{
T(\alpha_0) =
\begin{cases}
  {\gamma Ca_\s^{\half}}\cdot\f{M\delta\log(n)}{\sqrt{n}}\ \ &\text{ if } p=\half,\\ \\
  \f{\gamma a_\s^p}{(1-2p)} \cdot \f{M\delta}{\sqrt{n}} - \f{Kpa_\s^{\half}}{(1-2p)} \cdot \f{M\delta}{n^{{1}/{4p}}}\ \ &\text{ if } 0<p\neq\half<1,
\end{cases}
\label{eq:cases}
}
for some fixed constant $C > 3 - \log(9a)$. Observe that when $0<p<\half$, the last term of Eq.~\eqref{eq:cases} is negative, and similarly when $\half<p<1$, the first term is negative. From this, we have that $T(\alpha_0) \le M\delta \xi(n,p)$ where
\eq{
\xi(n,p) =
\begin{cases}
  \f{\gamma a_\s^p}{(1-2p)} \cdot \f{1}{\sqrt{n}} & \text{ if } 0 < p < \half,\nn\\\nn\\
  {\gamma Ca_\s^{\half}}\cdot\f{\log(n)}{\sqrt{n}} & \text{ if }  p = \half,\nn\\ \nn\\
  \f{\gamma pa_\s^{\half}}{2p-1} \cdot \f{1}{n^{{1}/{4p}}} & \text{ if } \half < p < 1.
\end{cases}
}

Plugging this into Eq.~\eqref{eq:mcdiarmid1}, we have that with probability greater than $1-e^{-t}$,
\eq{
\supgdelta < M\delta\xi(n,p) + M\delta\sqrt{\f{2t}{n}}.
\label{eq:temp1}
}

From Eq.~\eqref{eq:pr-gdel}, this implies that
\eq{
\pr^{\otimes n}\pb{\Xv_{1:n} : \normh{\fn-\fo} > \delta} \le \pr^{\otimes n}\pb{\Xv_{1:n} :   \sup_{g \in \gdelta}\Delta_n(g) \ge \f{\mu\delta^2}{2}}.\nn
}
Thus, in Eq.~\eqref{eq:temp1}, by letting
\eq{
\f{\mu\delta^2}{2} = \pa{M\delta\xi(n,p) + M\delta\sqrt{\f{2t}{n}}},\nn
}
we have that with probability greater than $1-e^{-t}$,
\eq{
\normh{\fn-\fo} \le \f{2M}{\mu}\pa{\xi(n,p) + \sqrt{\f{2t}{n}}}.\nn
}
Observe that $\norminf{\fn-\fo} \le \kinf^\half\normh{\fn-\fo}$. For $0 < \alpha < 1$, by choosing $\delta_n$ as
\eq{
\delta_n = \f{2M\kinf^\half}{\mu}\pa{\xi(n,p) + \sqrt{\f{2\log(1/\alpha)}{n}}},\nn
}
we have that
\eq{
\pr^{\otimes n}\pb{\Xv_{1:n} : \norminf{\fn-\fo} \le \delta_n} > 1-\alpha.\nn
}
From the stability of persistence diagrams in Proposition \ref{prop:stability}, this implies that
\eq{
\pr^{\otimes n}\bigg\{\Xv_{1:n} : \Winf \Big( \dgm\left(\fn\right),\dgm\left(\fo\right) \Big) > \delta_n\bigg\} \le \alpha,\nn
}
yielding the desired result. \qed


\section{Supplementary Results}
\label{supplementary-results}

In this section, we establish some results which play a key role in the proofs presented in Section \ref{proofs}.

\subsection{Properties of the Risk Functional $\J(g)$}

We establish some important properties of the risk functional, given by
\eq{
\J(g) = \int_{\R^d}\ell_g(\xv) \ d\pr(\xv) = \int_{\R^d}{\rho\pa{\normh{\phigox}}}  \ d\pr(\xv).\nn
}

The following result establishes that some important properties of the robust loss $\rho$ carry forward to $\J(g)$. (i) The Lipschitz property of $\rho$ is inherited by $\J(g)$, (ii) the convexity of $\rho$ is strengthened to guarantee that $\J(g)$ is strictly convex, and (iii) $\J(g)$ is strongly convex with respect to the $\normh{\cdot}$--norm around its minimizer.


\begin{proposition}[Convexity and Lipchitz properties of $\J$]
  Under assumptions $\pa{\mathcal{A}1}-\pa{\mathcal{A}3}$,
  \begin{enumerate}[label=(\roman*)]
    \item{The risk functionals $\J(g)$ and $\J_n(g)$ are $M$-Lipschitz w.r.t. $\norm{\cdot}_\HH$.}

    \item Furthermore, if $\rho$ is convex, $\J(g)$ and $\J_n(g)$ are strictly convex.
    \item Additionally, under assumption $\pa{\mathcal{A}4}$, for $\fo = \arginf_{g\in\G}\J(g)$, the risk functional \mbox{satisfies} the strong convexity condition
    \eq{
    \J(g) - \J(\fo) \ge \f{\mu}{2}\norm{\fo-g}^2_\HH,\nn
    }
    for $\mu = 2\min\pb{\varphi\pa{2\kinf^\half}, \rho''\pa{2\kinf^\half}}$.
  \end{enumerate}
  \label{lemma:lipschitz}
\end{proposition}

\prf{
\textbf{Lipschitz property.} Observe that,
\eq{
\abs{\ell_{g_1}(\xv) - \ell_{g_2}(\xv)} &= \abs{\rho\pa{\norm{\Phis(\xv)-g_1}_\HH} - \rho\pa{\norm{\Phis(\xv)-g_2}_\HH}} \nn\\
&\le M \abs{\norm{\Phis(\xv)-g_1}_\HH - \norm{\Phis(\xv)-g_2}_\HH}\nn\\
&\le M \norm{g_1-g_2}_\HH \ ,\nn
}
where the first inequality follows from the fact that $\rho$ is $M$-Lipschitz and the last inequality follows from reverse triangle inequality. This shows that the loss functions $\ell_g(\cdot)$ are
$M$-Lipschitz with respect to $g$. For the risk functionals, we have that,
\eq{
\abs{\J(g_1) - \J(g_2)} &= \abs{\int\limits_{\R^d}{\pa{\ell_{g_1}(\xv) - \ell_{g_2}(\xv)} d\pr(\xv)}}\nn\\
&\le {\int\limits_{\R^d}{\Abs{\ell_{g_1}(\xv) - \ell_{g_2}(\xv)} d\pr(\xv)}}\nn\\
&\le M \norm{g_1-g_2}_\HH,\nn
}
where the first inequality follows from Jensen's inequality. This verifies that $\J(g)$ is $M$-Lipchitz. The proof for $\J_n(g)$ is identical.

\textbf{Strict Convexity.} We begin by establishing that for translation invariant kernels $\normh{\Phis(\xv)-\cdot}$ is strictly convex. Suppose $g_1,g_2 \in \HH \cap \D$ and $\lambda \in (0,1)$, and let $g = (1-\lambda)g_1 + \lambda g_2$. Then
\eq{
\normh{\Phis(\xv) - g}^2 &= \normh{(1-\lambda)(\Phis(\xv) - g_1) + \lambda(\Phis(\xv) - g_2)}^2\nn\\
&= (1-\lambda)^2\normh{\Phis(\xv)-g_1}^2\nn\\
& \ \ \ \ + \lambda^2\normh{\Phis(\xv)-g_2}^2 + 2\lambda(1-\lambda)\Innerh{\Phis(\xv)-g_1,\Phis(\xv)-g_2}.
\label{eq:rkhs-convexity}
}

From Cauchy-Schwarz inequality, we know that
\eq{
\Innerh{\Phis(\xv)-g_1,\Phis(\xv)-g_2} \le \normh{\Phis(\xv)-g_1}\normh{\Phis(\xv)-g_2}.\nn
}

In the following, we argue that for translation invariant kernels,
\eq{
\Innerh{\Phis(\xv)-g_1,\Phis(\xv)-g_2} < \normh{\Phis(\xv)-g_1}\normh{\Phis(\xv)-g_2},
\label{eq:cauchy-schwarz}
}
for $g_1\neq g_2$. On the contrary, suppose
\eq{
\Innerh{\Phis(\xv)-g_1,\Phis(\xv)-g_2} = \normh{\Phis(\xv)-g_1}\normh{\Phis(\xv)-g_2}\nn
}
holds. Then this implies that there is a function $a(\xv)$, depending only on $g_1$ and $g_2$, such that $a(\xv) \neq 0$ for $\xv \in \R^d$ and
\eq{
\Phis(\xv)-g_1 = a(\xv) \pa{\Phis(\xv) - g_2}.\nn
}
Rearranging the terms this implies that
\eq{
\Phis(\xv) = \f{g_1 - a(\xv)g_2}{1-a(\xv)} = (1+b(\xv))g_1 + b(\xv)g_2,\nn
}
where $b(\xv) = -a(\xv)/(1-a(\xv))$ also does not vanish on $\xv \in \R^d$. For $\xv,\yv \in \R^d$, from the reproducing property we have
\eq{
\K(\xv,\yv) &= \Innerh{\Phis(\xv),\Phis(\yv)}\nn\\
&= \Innerh{g_1+b(\xv)(g_1+g_2),g_1+b(\yv)(g_1+g_2)}\nn\\
&= b(\xv)b(\yv)\normh{g_1+g_2}^2 + \pa{b(\xv)+b(\yv)}\innerh{g_1,g_1+g_2} + \normh{g_1}^2.\nn
}
Note that because the kernel is translation invariant, i.e., $\K(\xv,\xv) = \K(\yv,\yv) = \s^{-d}\psi(0)$, this must imply that
\eq{
0&=\pa{b(\xv)^2 - b(\yv)^2}\normh{g_1+g_2}^2 + 2(b(\xv)-b(\yv))\innerh{g_1,g_1+g_2}\nn\\
&=\pa{b(\xv)-b(\yv)}\pa{(b(\xv)+b(\yv))\normh{g_1+g_2}^2 + 2\innerh{g_1,g_1+g_2}}.\nn
}
Since $b(\xv)$ and $b(\yv)$ are nonvanishing, the above equation is satisfied only when $b(\xv) = b(\yv)$. This implies that $\K(\xv,\yv)$ is constant for all $\yv$, giving us a contradiction. Thus, we have that Eq.~\eqref{eq:cauchy-schwarz} holds. Plugging this back in Eq.~\eqref{eq:rkhs-convexity} we get that for $\lambda \in (0,1)$ and $g = (1-\lambda)g_1 + \lambda g_2$,
\eq{
\normh{\Phis(\xv)-g} < (1-\lambda)\normh{\Phis(\xv)-g_1} + \lambda \normh{\Phis(\xv)-g_2}.\nn
}
Since, $\rho$ is strictly increasing and convex, this implies that
\eq{
\ell_g(\xv) < (1-\lambda)\ell_{g_1}(\xv) + \lambda \ell_{g_2}(\xv).\nn
}
The map $\ell_g(\cdot) \mapsto \pr \ell_g$ is a linear operator, and $\ell_g$ is strictly convex in $g$, this implies that $\J(g)$ is also strictly convex in $g$. The same holds for $\J_n(g)$.

\textbf{Strong Convexity around the minimizer.} We now turn our attention to the strong convexity property. For this, we first show that $\J(g)$ is twice \gat{} differentiable. Let $g,h\in \G$, then the second \gat{} derivative of the loss $\ell_g(\xv) = \rho\pa{\norm{\Phis(\xv)-g}_\HH}$ at $g$ in the direction $h$ is given by,
\eq{
\delta^2\ell(\xv,g;h) &= \f{d^2}{d\alpha^2}\ell(\xv,g+\alpha h)\Big\vert_{\alpha=0}\nn\\
&= \f{d^2}{d\alpha^2}\rho\pa{\norm{\Phis(\xv)-g -\alpha h}_\HH}\Big\vert_{\alpha=0}\nn\\
&= \f{d}{d\alpha} \pc{\varphi\pa{\norm{\Phis(\xv)-g-\alpha h}_\HH}\pa{-\innerh{\Phis(\xv)-g,h} + \alpha \normh{h}^2}}\Big\vert_{\alpha=0}\nn\\
&= \varphi\pa{\normh{\Phis(\xv)-g}}\normh{h}^2 + \innerh{\Phis(\xv)-g,h}^2 \f{\varphi'\pa{\normh{\Phis(\xv)-g}}}{\normh{\Phis(\xv)-g}}\nn\\
&= \varphi\pa{z(\xv,g)}\normh{h}^2 +  \normh{h}^2 \lambda(\xv,g,h) z(\xv,g) \varphi'\pa{z(\xv,g)},
\label{eq:second-deriv}
}
where for a fixed $g \in \G$, in the interest of brevity, we define $z(\xv,g) = \normh{\Phis(\xv)-g}$ and
\eq{
\lambda(\xv,g,h) = \Inner{\f{\Phis(\xv)-g}{\normh{\Phis(\xv)-g}},\f{h}{\normh{h}}}_\HH^2 \in [0,1].\nn
}
Observe that $z\varphi'(z) = \rho''(z) - \varphi(z)$, thus Eq.~\eqref{eq:second-deriv} becomes
\eq{
\delta^2\ell(\xv,g;h) &= \normh{h}^2 \left(\pa{1-\lambda(\xv,g,h)}\varphi\pa{z(\xv,g)} + \lambda(\xv,g,h) \rho''\pa{z(\xv,g)}\right).\nn
}

From assumption $\pa{\mathcal{A}4}$ we have that $\rho''$ and $\varphi$ are nonincreasing, and
\eq{
z(\xv,g) = \normh{\Phis(\xv)-g} \le 2\kinf^\half.\nn
}

Thus, we have that
\eq{
\delta^2\ell(\xv,g;h) \ge c \normh{h}^2,
\label{eq:gateaux1}
}
where
\eq{
c = \min\pb{\varphi\pa{2\kinf^\half},\rho''\pa{2\kinf^\half}}.\nn
}
We also note that $\delta^2\ell(\xv,g;h)$ is bounded above. To see this, note that from assumption $\pa{\mathcal{A}4}$, $\rho''$ and $\varphi$ are bounded and nonincreasing. Consequently, for $\lambda(\xv,g,h) \in (0,1)$ and
\eq{
C = \max\pb{\rho''(0),\varphi(0)} < \infty,\nn
}
from Eq.~\eqref{eq:second-deriv} we have that
\eq{
\delta^2\ell(\xv,g;h) \le C \normh{h}^2 < \infty.\nn
}

The \gat{} derivative of $\J(g)$ is, then, given by
\eq{
\delta^2\J(g;h) &= \f{d^2}{d\alpha^2}\J(g+\alpha h)\Big\vert_{\alpha=0} = \f{d^2}{d\alpha^2}\int\limits_{\R^d}{\ell(\xv,g+\alpha h) \ d\pr(\xv)}\Big\vert_{\alpha=0}\nn\\
&= \int\limits_{\R^d}{\f{d^2}{d\alpha^2}\ell(\xv,g+\alpha h) \ d\pr(\xv)}\Big\vert_{\alpha=0}\nn\\
&= \int\limits_{\R^d}{\delta^2\ell(\xv,g;h) \ d\pr(\xv)}.\nn
}
The exchange of the derivative and integral in the second line follows from the dominated convergence theorem since $\abs{\delta^2\ell(\xv,g;h)}$ is bounded. This confirms the \gat{} differentiability of $\J(g)$. From Eq.~\eqref{eq:gateaux1} we have
\eq{
\delta^2\J(g;h) = \int\limits_{\R^d}{\delta^2\ell(\xv,g;h) \ d\pr(\xv)} \ge c \normh{h}^2.
\label{eq:gateaux2}
}
For $\fo = \arginf_{g\in\G}\J(g)$ and $g\in \G$, we proceed to show the strong-convexity guarantee. Let $h = g - \fo$. From the first-order Taylor approximation for $\J(g)$ we have,
\eq{
\J(g) = \J(\fo) + \delta\J(\fo,h) + R_2(\fo,h),\nn
}
where the first \gat{} derivative, $\delta\J(\fo,h) = 0$ for all $h$ since $\fo$ is the unique minimizer of $\J(g)$ and the remainder term $R_2(\fo,h)$ is given by
\eq{
R_2(\fo,h) &= \f{1}{2}\int_0^1{(1-t)\delta^2\J(\fo + th;h) \ dt}\nn\\
&\ge \f{c}{2} \normh{h}^2 \int_0^1{(1-t)dt} = \f{c}{4}\normh{h}^2,\nn
}
where the inequality follows from Eq.~\eqref{eq:gateaux2}. As a result, for any $g \in \G$ and  $\mu = \f{c}{2}$ we have that
\eq{
\J(g) - \J(\fo) \ge \f{\mu}{2}\normh{g-\fo}^2,\nn
}
yielding the desired result.
}

We now turn to examining the behaviour of the risk functional $\J(g)$ w.r.t. the underlying probability measure $\pr$. For $0 \le \e \le 1$ and $\xv \in \R^d$, let $\prx = (1-\e)\pr + \e\delta_{\xv}$ be a perturbation curve, as defined in Theorem \ref{thm:influence}. The risk functional associated with $\prx$ is given by
\eq{
\Jx(g) = \prx\ell_g = (1-\e)\J(g) + \e\rho\pa{\normh{\phigox}},\nn
}
and $\fox = \inf_{g\in \G}\Jx(g)$ is the minimizer. The convergence of $\fox$ to $\fo$ can be studied by examining the convergence of $\Jx$ to $\J$. Specifically, under conditions on $\J$ and $\Jx$, it can be shown that $\normh{\fox-\fo} \rightarrow 0$ as $\e \rightarrow 0$. The machinery we use here uses the notion of $\Gamma$--convergence, which is defined as follows.

\begin{definition}[$\Gamma$ convergence]
  Given a functional $F: \X \rightarrow \R \cup \pb{\pm \infty}$ and a sequence of functionals $\pb{F_n}_{n\in \N}$, $F_n \stackrel{\Gamma}{\rightarrow} F$ as $n\rightarrow \infty$ when
\begin{enumerate}[label=(\roman*)]
  \item $F(\xv) \le \liminf\limits_{n\rightarrow \infty} F_n(\xv_n)$ for all $\xv \in \X$ and every $\{\xv_n\}_{n\in\N}$ such that ${d(\xv_n,\xv) \rightarrow 0}$;
  \item For every $\xv \in \X$, there exists $\{\xv_n\}_{n \in \N}$, $d(\xv_n,\xv) \rightarrow 0$ such that $F(\xv) \ge \limsup\limits_{n\rightarrow \infty} F_n(\xv_n)$.
\end{enumerate}
\end{definition}

The following result shows that the sequence of functionals $\pb{\Jx}$ $\Gamma$--converges to $\J$.

\begin{proposition}[$\Gamma$--convergence of $\Jx$ to $\J$]
  Under assumptions $\pa{\mathcal{A}1}$--$\pa{\mathcal{A}3}$,
  \eq{
  {\Jx(g) \stackrel{\Gamma}{\rightarrow} \J(g)} \ \ \ \text{ as } \ \e \rightarrow 0.\nn
  }
  \label{lemma:gamma}
\end{proposition}
\prf{
  Let $g \in \G$ and $\pb{g_\e}_{\e > 0}$ be a sequence in $\G$ such that $\normh{g_\e-g} \rightarrow 0$ as $\e \rightarrow 0$. In order to verify $\Gamma$--convergence we first show that the following holds
  \eq{
  \lim_{\e \rightarrow 0}\Abs{\Jx(g_\epsilon) - \J(g)} = 0.\nn
  }
  For $\e > 0$, using the triangle inequality we have that
  \eq{
  \Abs{\Jx(g_\epsilon) - \J(g)} &\le \Abs{\J(g_\e) - J(g)} + \Abs{\Jx(g_\e)-\J(g_\e)}\nn\\
  &\stackrel{(i)}{\le} M\normh{g_\e-g} + \Abs{\Jx(g_\e)-\J(g_\e)}\nn\\
  &\stackrel{(ii)}{\le} M\normh{g_\e-g}  + \epsilon\cdot\Abs{\J(g) - \rho\pa{\normh{\phigox}}},\nn
  }
  where (i) uses the fact that $\J(g)$ is $M$--Lipschitz from Proposition \ref{lemma:lipschitz}, and (ii) uses the fact that
  \eq{
  \Jx(g) = (1-\e)\J(g) + \e \rho\pa{\normh{\phigox}}.\nn
  }
  Since $\normh{g_\e-g} \rightarrow 0$ as $\e \rightarrow 0$ we have
  \eq{
  \lim_{\e \rightarrow 0}\Abs{\Jx(g_\epsilon) - \J(g)} \le M\lim_{\e \rightarrow 0}\normh{g_\e-g} + \lim_{\e \rightarrow 0}\epsilon\cdot\Abs{\J(g) - \rho\pa{\normh{\phigox}}} = 0.\nn
  }
  Since $\Jx$ and $\J$ are continuous, using \cite[Remark 4.8]{dal2012introduction} it follows that $\Jx(g) \stackrel{\Gamma}{\rightarrow} \J(g)$.
}

Now, we examine the coercivity of the sequence $\pb{\Jx}$.

\begin{definition}[Equi-coercivity]
  A sequence of functionals $\pb{F_n}_{n\in \N} : \X \rightarrow \R \cup \pb{\pm \infty}$ is~said to be equi-coercive if for every $t \in \R$, there exists a compact set $K_t \subseteq \X$ such that ${\pb{\xv \in \X : F_n \le t} \subseteq K_t}$ for every $n \in \N$.
\end{definition}

The following result shows that the sequence $\pb{\Jx}$ is equi-coercive.

\begin{proposition}[Equi-coercivity of $\Jx$]
    Under assumptions $\pa{\mathcal{A}1}$--$\pa{\mathcal{A}3}$, the sequence of \mbox{functionals} $\pb{\Jx}$ is equi-coercive.
    \label{lemma:equi}
\end{proposition}

\prf{
 For $0 < \epsilon < 1$, $\xv \in \R^d$ and $g \in \G$, we have that
 \eq{
    \Jx(g) = (1-\epsilon)\J(g) + \epsilon \rho\pa{\normh{\phigox}}.\nn
 }
 From \cite[Proposition 7.7]{dal2012introduction} in order to show that the sequence of functionals $\pb{\Jx}$ is equi-coercive, it suffices to show that there exists a lower semicontinuous, coercive functional $F:\HH \rightarrow \R \cup \pb{\pm \infty}$ such that $F \le \Jx$ for every $\epsilon \ge 0$. To this end consider the functional
 \eq{
 F(g) = \min\pb{\J(g),\rho\pa{\normh{\phigox}}}.\nn
 }
 As $\Jx$ is a convex combination of $\J(\cdot)$ and $\rho\pa{\normh{\Phis(\xv) - \cdot}}$, it implies that $F \le \Jx$ for every $\e \ge 0$. Additionally, because $\J(\cdot)$ and $\rho\pa{\normh{\Phis(\xv) - \cdot}}$ are both continuous, it follows that $F$ is also continuous, and, therefore, lower semicontinuous.

 We now verify that $F$ is coercive. Since $\rho$ is strictly increasing we have that
 \eq{
  \rho\pa{\normh{\phigox}} \rightarrow \infty  \ \ \ \ \text{ as } \ \normh{g} \rightarrow \infty,\nn
 }
 verifying that $\rho\pa{\normh{\Phis(\xv) - \cdot}}$ is coercive. Next, from the reverse triangle inequality we have that
 \eq{
 \normh{\phigox} &\ge \Abs{\normh{\Phis(\xv)} - \normh{g}} = \Abs{\sqrt{\K(\xv,\xv)} - \normh{g}}.\nn
 }
 Observe that $\K(\xv,\xv) = \kinf$, and because $\rho$ is strictly increasing we have
 \eq{
    \rho\pa{\Abs{\kinf^\half - \normh{g}}} \le \rho\pa{\normh{\phigox}}.\nn
 }
 Taking expectations on both sides w.r.t. $\pr$,
 \eq{
 \rho\pa{\Abs{\kinf^\half - \normh{g}}} \le \int_{\R^d}\rho\pa{\normh{\phigox}}d\pr(\xv) = \J(g).\nn
 }
Since
\eq{
\rho\pa{\Abs{\kinf^\half - \normh{g}}} \rightarrow \infty  \ \ \ \ \text{ as } \ \normh{g} \rightarrow \infty,\nn
}
it implies that $\J(g)$ is coercive as well. It follows from this that $F$ is coercive, and the sequence of functionals $\pb{\Jx}$ is equi-coercive.
}

Propositions \ref{lemma:gamma} and \ref{lemma:equi} together imply, from the fundamental theorem of $\Gamma$-convergence \cite[Theorem 7.8]{dal2012introduction}, that the sequence of minimizers associated with $\pb{\Jx}$ converge to the minimizer of $\J$, i.e.,
\eq{
{\normh{\fox-\fo} \rightarrow 0} \ \ \  \ \text{ as }  \ \e \rightarrow 0.\nn
}

\subsection{Some Additional Results}

Next, we note an important property of the hypothesis class, $\G = \HH \cap \D$. The elements of $\G$ can be shown to have their $\norm{\cdot}_\infty$--norm related their  $\norm{\cdot}_\HH$--norm.

\begin{lemma}[{\cite[Lemma 6]{vandermeulen2013consistency}} and {\cite[Proposition 5.1]{sriperumbudur2016optimal}}]
  For every $g \in \HH \cap \D$,
  \eq{
  \norm{g}^2_\HH \le \norm{g}_\infty \le \norm{\K}^{\half}_\infty \norm{g}_\HH.\nn
  }
  \label{lemma:supnorm}
\end{lemma}

\renewcommand{\Qs}{R_\s}

The following result, which is essentially the population analogue of \cite[Lemma 7]{vandermeulen2013consistency}, guarantees that for small enough $\s > 0$, there exists $0< \delta < 1$ such that $\fo$ is contained in the RKHS ball $B_{\HH}(\zerov,\delta\nus)$, where for brevity we denote $\nus = \kinf^{1/2}$. We provide the proof for completeness, however, the proof uses exactly the same ideas from \cite{vandermeulen2013consistency}. For notational convenience, we also define $\psis(\norm{\xv-\yv}_2) = \K(\xv,\yv) = \s^{-d}\psi\pa{\norm{\xv-\yv}_2/\s}$.

\begin{lemma}
  Let $\pr \in \mathcal{M}(\R^d)$ and $\fo$ be the robust KDE for $\s > 0$. For sufficiently small $\s > 0$, there exists $0 < \delta < 1$ such that $\fo \in B(\zerov,\delta\nus)$.
  \label{lemma:vanishing}
\end{lemma}

\prf{
  For $\pr \in \M(\R^d)$, and $\G = \HH \cap \D$, consider the map $\Ts : \G \rightarrow \G$ given by
  \eq{
  \Ts(g) = \int_{\R^d}{\f{\varphi\pa{\normh{\phigox}}}{\int_{\R^d}\varphi\pa{\normh{\phigoy}} d\pr(\yv)} \K(\cdot,\xv) \ d\pr(\xv) } = \int_{\R^d}{\K(\cdot,\xv)\ws(\xv)d\pr(\xv)},\nn
  }
  for each $g \in \G$. Observe that $\ws \in L_1(\pr)$ is a non-negative function such that
  \eq{
  \int_{\R^d}{\ws(\xv)d\pr(\xv)} = 1.
  \label{eq:ws1}
  }
  Let $S_\s = \textup{Im}(\Ts) \subset \G$. It follows from \cite[Page 11]{vandermeulen2013consistency} that the robust KDE, $\fo = \arginf_{g \in \G}\J(g)$, is the fixed point of the map $\Ts$ and therefore $\fo \in S_\s$. For a small $\e > 0$, from \cite[Lemma 12; Corollary 13]{vandermeulen2013consistency} there exist $r,s > 0$ such that $\pr(B(\xv,r)) \le \e$ and $\pr(B(\xv,r+s)\setminus B(\xv,r)) \le \e$ for all $\xv \in \R^d$. This implies that $\pr(B(\xv,r+s)^c) > 1-2\e$. We point out that the constant $\e$ chosen here is related to $\sqrt{9/10}$ used by \cite{vandermeulen2013consistency} as $\sqrt{1-\e} = \sqrt{9/10}$, which, as remarked by the authors, was chosen simply for convenience. Define the sets $B_\s = B_{\HH}(\zerov,\nus\sqrt{1-\e})$, and let
  \eq{
  \Qs \defeq S_\s \cap B_\s^c.\nn
  }

  In what follows we will show that $\fo$ does not lie in $\Qs$. To this end, let $g = \arginf_{h \in \Qs}\J(h)$. It suffices to show that $\J(g) > \J(\zerov) > \J(\fo)$. Since $g \in \Qs$, it must follow that
  \eq{
  (1-\e)\nus^2 < \normh{g}^2 \le \norminf{g} = g(\zv),
  \label{eq:g-ineq}
  }
  for some $\zv \in \R^d$, where the second inequality follows from Lemma \ref{lemma:supnorm}. Since $g \in S_\s$, there exists a non-negative function $\ws$ satisfying Eq.~\eqref{eq:ws1}, such that $g = \int_{\R^d}{\ws(\xv)\K(\cdot,\xv)d\pr(\xv)}$. Therefore,
  \eq{
  (1-\e)\nus^2 \le g(\zv) &= \int_{\R^d}\K(\zv,\xv)\ws(\xv)d\pr(\xv)\nn\\
  &= \int_{B(\zv,r)}{\K(\zv,\xv)\ws(\xv)d\pr(\xv)} + \int_{B(\zv,r)^c}{\K(\zv,\xv)\ws(\xv)d\pr(\xv)}\nn\\
  &\stackrel{(i)}{\le} \nus^2 \int_{B(\zv,r)}{\ws(\xv)d\pr(\xv)} +  \psis(r) \underbrace{\int_{B(\zv,r)^c}{\ws(\xv)d\pr(\xv)}}_{\le 1}\nn\\
  &\stackrel{(ii)}{\le} \nus^2 \int_{B(\zv,r)}{\ws(\xv)d\pr(\xv)} +  \psis(r),
  \label{eq:psi-bound1}
  }
  where (i) follows from the fact that $\sup_{B(\zv,r)^c}\K(\zv,\xv) = \psis(r)$ and (ii) follows from Eq.~\eqref{eq:ws1}. From \cite[Lemma 7]{vandermeulen2013consistency}, there exists $\s$ small enough such that $ \psis(r) < \f{\e}{2}\nus^2$. Plugging this back in Eq.~\eqref{eq:psi-bound1} we get
  \eq{
  \int_{B(\zv,r)}{\ws(\xv)d\pr(\xv)} \ge \pa{1-\f{3\e}{2}}.
  \label{eq:term1}
  }
  Additionally,
  \eq{
  \sup\limits_{\yv \in B(\zv,r+s)^c}g(\yv) &= \sup\limits_{\yv \in B(\zv,r+s)^c}\Bigl( \int\limits_{B(\zv,r)}{\K(\yv,\xv)\ws(\xv)d\pr(\xv)} + \int\limits_{B(\zv,r)^c}{\K(\yv,\xv)\ws(\xv)d\pr(\xv)} \Bigr)\nn\\
  &\le \sup\limits_{\yv \in B(\zv,r+s)^c}\sup\limits_{\xv \in B(\zv,r)}\K(\yv,\xv) \int\limits_{B(\zv,r)}{\ws(\xv)d\pr(\xv)} \nn\\
  & \ \ \ \  \ \ \ + \sup\limits_{\yv \in B(\zv,r+s)^c}\sup\limits_{\xv \in B(\zv,r)}\K(\yv,\xv)\int\limits_{B(\zv,r)^c}{\ws(\xv)d\pr(\xv)}\nn\\
  &\le  \psis(s) + \nus^2 \int\limits_{B(\zv,r)^c}\ws(\xv)d\pr(\xv).\nn
  }
  For a choice of $\tau > 0$, there is $\s$ small enough satisfying $\psis(s) \le \tau$ such that from Eq.~\eqref{eq:term1}
  \eq{
  \sup\limits_{\yv \in B(\zv,r+s)^c}g(\yv) \le \tau + {\f{3\e}{2}}\nus^2.
  \label{eq:term2}
  }
  Then we have that
  \eq{
    \J(g)  &= \int\limits_{\R^d}{\rho\pa{\normh{\phigox}} d\pr(\xv)}\nn\\
    &= \int\limits_{B(\zv,r+s)}{\rho\pa{\normh{\phigox}} d\pr(\xv)} + \int\limits_{B(\zv,r+s)^c}{\rho\pa{\normh{\phigox}} d\pr(\xv)}\nn\\
    &{\ge} \int\limits_{B(\zv,r+s)^c}{\rho\pa{\normh{\phigox}} d\pr(\xv)}\nn\\
    &=\int\limits_{B(\zv,r+s)^c}{\rho\pa{\sqrt{\nus^2+ \normh{g}^2 - 2 \innerh{g,\Phis(\xv)}}} d\pr(\xv)}\nn\\
    &\ge \int\limits_{B(\zv,r+s)^c}{\rho\pa{\sqrt{\nus^2+ \normh{g}^2 - 2 \sup\limits_{\yv \in B(\zv,r+s)^c}g(\yv) } } d\pr(\xv)}.\nn
  }
  Plugging in Equations \eqref{eq:term2} and \eqref{eq:g-ineq} we get
  \eq{
  \J(g) \ge (1-2\e)\rho\pa{\sqrt{(2-4\e)\nus^2 - 2\tau}}.\nn
  }
  Since $\rho$ is assumed to be strictly convex, this implies that $\rho'$ is strictly increasing. Additionally, from $\pa{\mathcal{A}2}$ we have that $\rho'$ is bounded. This implies that, for any $0 < \alpha < \norminf{\rho'}$, there is $\beta > 0$ such that $\rho'(z) > \norminf{\rho'} - \alpha$ for all $z > \beta$.  Using \cite[Eq.~(11)]{vandermeulen2013consistency}, we have
  \eq{
  \rho\pa{\sqrt{(2-4\e)\nus^2 - 2\tau}} &= \int_{0}^{(2-4\e)\nus^2 - 2\tau}{\rho'(z)dz}\nn\\
  &\ge \int_{\beta}^{(2-4\e)\nus^2 - 2\tau}{\rho'(z)dz}\nn\\
  &\ge \int_{\beta}^{(2-4\e)\nus^2 - 2\tau}{\pa{\norminf{\rho'} - \alpha} dz}\nn\\
  &\ge \pa{\norminf{\rho'}-\alpha}\pa{\sqrt{(2-4\e)\nus^2 - 2\tau} - \beta}.\nn
  }
  Without loss of generality, we can assume $\norminf{\rho'}=1$. Choosing $\alpha$, $\tau$ and $\sigma$ small enough we obtain
  \eq{
  \J(g) \ge \nus.\nn
  }
  Now we note that
  \eq{
    \J(\zerov) &= \int_{\R^d}\rho\pa{\normh{\Phis(\xv)}}d\pr(\xv)\nn\\
    &= \rho\pa{\nus}\nn\\
    &= \rho(0) + \int_{0}^{\nus}\rho'(z)dz\nn\\
    & \ {\le}  \  \rho(0) + \norminf{\rho'}\int_0^{\nus}{dz} = \nus.\nn
  }
  Thus, we obtain that $\J(g) > \J(\zerov)$. We have $g = \arginf_{h\in\Qs}\J(h)$ and $\fo = \arginf_{h \in \G}\J(h)$, and, additionally we know that $\fo \neq \zerov$. It follows that since $\J(\fo) \le \J(\zerov) < \J(g)$, then $\fo \notin \Qs$ as $\s \rightarrow 0$. Taking $\delta = \sqrt{1-\e}$, we get the desired result.
}

\renewcommand{\Qs}{Q_\s}


\section{Supplementary Results for the Persistence Influence}
\label{influence}

In this section, we collect the proofs for the results on persistence influence established in Section~\ref{robustness}. The following result shows that when $\varphi$ is nonincreasing, the persistence influence in Eq.~\eqref{eq:influence-rkde} can be written in a more succinct form.

\begin{proposition}
  Under the conditions of Theorem \ref{thm:influence}, if $\varphi$ is nonincreasing, then the persistence influence of $\xv \in \R^d$ on $\dgm\pa{\fo}$ satisfies
  \eq{
  \Psi\pa{\fo;\xv} \le \kinf^\half \ws(\xv) \normh{\phifox},\nn
  }
  where $\ws$ is the measure of inlyingness from Eq.~\eqref{eq:weight}.
\end{proposition}

\prf{
From Theorem~\ref{thm:influence} we have that the persistence influence satisfies
\eq{
\Psi\pa{\fo;\xv} \le \kinf^\half  \rho'\pa{\normh{\phifox}}\left(\int_{\R^d}\zeta{\pa{\normh{\phifo}}}d\pr(\yv)\right)^{-1},
\label{eq:persinf}
}
where $\zeta(z) = \varphi(z) - z\varphi'(z)$. When $\varphi$ is nonincreasing, observe that $z\varphi'(z) \le 0$ for all $0 \le z < \infty$. Consequently, $\zeta$ can be bounded below by $\varphi$, and the r.h.s.~in Eq.~\eqref{eq:persinf} can be bounded above by
\eq{
\Psi\pa{\fo;\xv}
&\stackrel{(i)}{\le} \kinf^\half \f{ \rho'\pa{\normh{\phifox}}}{\int_{\R^d}\varphi{\pa{\normh{\phifo}}}d\pr(\yv)}\nn\\
&\stackrel{(ii)}{=} \kinf^\half \f{ \varphi\pa{\normh{\phifox}}}{\int_{\R^d}\varphi{\pa{\normh{\phifo}}}d\pr(\yv)} \normh{\phifox}\nn\\
&\stackrel{(iii)}{=} \kinf^\half \ws(\xv) \normh{\phifox},\nn
}
where (i) follows from the fact that $\zeta(z) \ge \varphi(z)$, (ii) follows from the definition of $\varphi$, i.e., ${\rho'(z) = z\varphi(z)}$, and (iii) follows from the definition of $\ws$ in Eq.~\eqref{eq:weight}, yielding the desired result.
}

The following result establishes the bound for the distance-to-measure described in Eq.~\eqref{eq:influence-dtm}.

\begin{proposition}
For $\pr \in \mathcal{M}(\R^d)$, the persistence influence for the distance-to-measure function is given by
\eq{
\Psi\pa{d_{\pr,m};\xv} \le \f{2}{\sqrt{m}}\sup\pb{\Abs{f(\xv) - \int_{\R^d}f(\yv)d\pr(\yv)} : \norm{\nabla f}_{L_2(\pr)} \le 1}\nn
}
where $\norm{\nabla f}_{L_2(\pr)}$ is a modified, weighted Sobolev norm \cite{villani2003topics,peyre2018comparison}.
\label{prop:dtm-influence}
\end{proposition}

\proof{
From \cite[Theorem 3.5]{chazal2011geometric} the following stability result holds:
\eq{
\norminf{d_{\pr,m} - d_{\prx,m}} \le \f{1}{\sqrt{m}}W_2\pa{\pr,\prx}.\nn
}
From \cite[Theorem 1]{peyre2018comparison} we have that
\eq{
W_2\pa{\pr,\prx} \le 2 \norm{\pr-\prx}_{\dot{H}\inv(\pr)},\nn
}
where the weighted, homogeneous Sobolev norm $\norm{\cdot}_{\dot{H}\inv(\mu)}$ for a signed measure $\nu$ w.r.t. a positive measure $\mu$ is given by
\eq{
\norm{\nu}_{\dot{H}\inv(\mu)} = \sup\pb{\Abs{\int_{\R^d}f(\xv)d\nu(\xv)} : \norm{\nabla f}_{L_2(\mu)} \le 1}.\nn
}
Observe that $\pr-\prx = \epsilon\pa{\delta_{\xv}-\pr}$ and since $\norm{\cdot}_{\dot{H}\inv(\mu)}$ defines a norm, we have that
\eq{
\lim_{\epsilon \rightarrow 0}\f{1}{\epsilon} \norminf{d_{\pr,m} - d_{\prx,m}} &\le \f{1}{\sqrt{m}} \lim_{\epsilon \rightarrow 0}\f{1}{\epsilon}W_2\pa{\pr,\prx}\nn\\
&\le  \f{2}{\sqrt{m}} \lim_{\epsilon \rightarrow 0} \f{1}{\epsilon}\norm{\epsilon\pa{\delta_{\xv}-\pr}}_{\dot{H}\inv(\pr)}\nn\\
&= \f{2}{\sqrt{m}} \norm{\pa{\delta_{\xv}-\pr}}_{\dot{H}\inv(\pr)}\nn\\
&= \f{2}{\sqrt{m}}\sup\pb{\Abs{f(\xv) - \int_{\R^d}f(\yv)d\pr(\yv)} : \norm{\nabla f}_{L_2(\pr)} \le 1}.\nn
}
From the stability for persistence diagrams, we have that
\eq{
\Psi\pa{d_{\pr,m};\xv} \le \lim_{\epsilon \rightarrow 0}\f{1}{\epsilon} \norminf{d_{\pr,m} - d_{\prx,m}}\nn
} and the result follows. \qed
}
\filbreak

\begin{figure}[H]
  \centering
  \begin{subfigure}[b]{0.32\linewidth}
    \includegraphics[width=\linewidth]{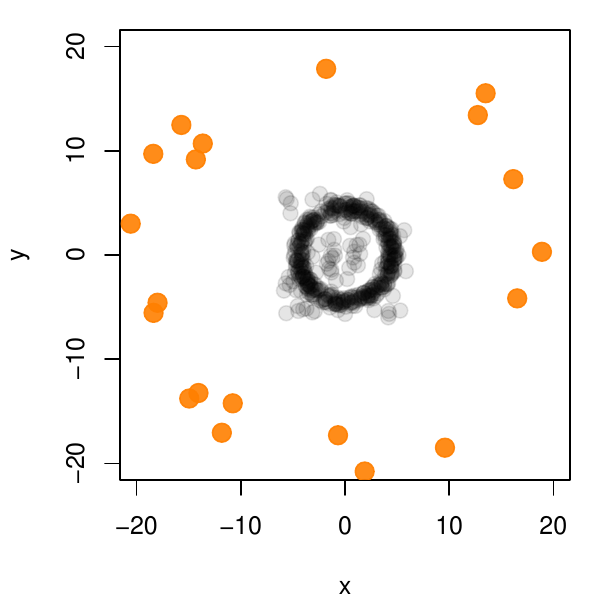}
    \caption{}
  \end{subfigure}
  \begin{subfigure}[b]{0.32\linewidth}
    \includegraphics[width=\linewidth]{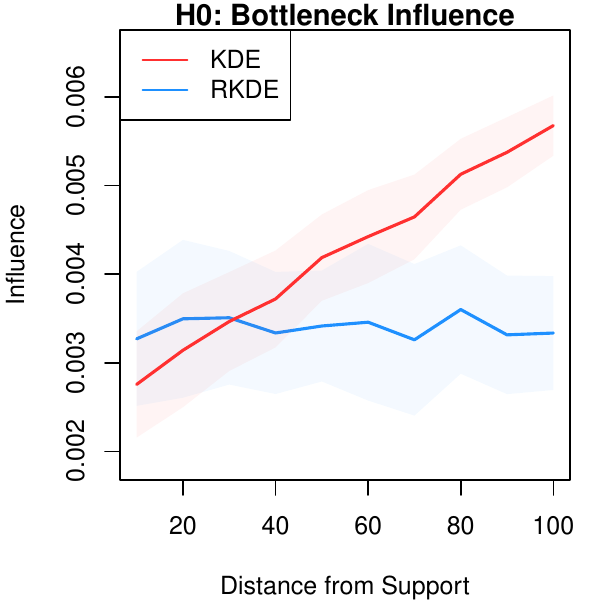}
    \caption{}
  \end{subfigure}
  \begin{subfigure}[b]{0.32\linewidth}
    \includegraphics[width=\linewidth]{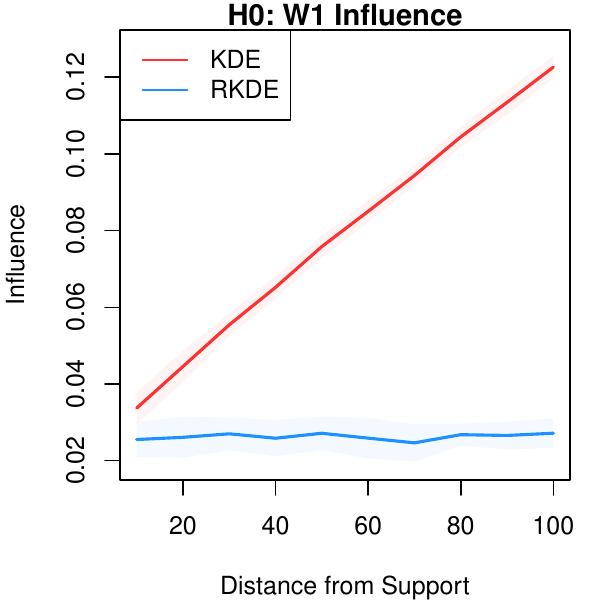}
    \caption{}
  \end{subfigure}
  \begin{subfigure}[b]{0.32\linewidth}
    \includegraphics[width=\linewidth]{./plots/influence/supnorm.pdf}
    \caption{}
  \end{subfigure}
  \begin{subfigure}[b]{0.32\linewidth}
    \includegraphics[width=\linewidth]{./plots/influence/h1b.pdf}
    \caption{}
  \end{subfigure}
  \begin{subfigure}[b]{0.32\linewidth}
    \includegraphics[width=\linewidth]{./plots/influence/h1w.pdf}
    \caption{}
  \end{subfigure}
  \caption{(a) An example of $\Xn$ in blue and the contamination $\Y_m$ when $r=10$. (d) The $L_{\infty}$ influence of $\Y_m$ on the KDE and robust KDE. (b, e) The bottleneck influence of $\Y_m$. (c, f) The $1$-Wasserstein influence of $\Y_m$ as the distance $r$ increases. }
  \label{fig:influence2}

\end{figure}

\textbf{Persistence-Influence Experiment} Points $\Xn$ are sampled from an annular region inside $\pc{-5,5}^2$ along with some uniform noise in the ambient space, corresponding to the black points in Figure~\ref{fig:influence2}~(a). $\Xn$ has interesting $1^{st}$-order homological features. We compute the robust KDE $\fn$ and the KDE $\barfn$ on the points $\Xn$ along with the corresponding persistence diagrams $\dgm\pa{\fn}$ and $\dgm\pa{\barfn}$. Outliers $\Y_m$ are added to the original points at a distance $r$ from the origin, the number of points roughly equal to $r$. Figure \ref{fig:influence2} (a) depicts these outliers in orange when $r=20$. The robust KDE $\fnm$ and $\barfnm$ are now computed on the composite sample $\Xn\cup\Y_m$ along with the persistence diagrams $\dgm\pa{\fnm}$ and $\dgm\pa{\barfnm}$. The bandwidth $\s(k)$ is chosen as the median distance to the $k^{th}$--nearest neighbour of each $\xv_i \in \Xn$, for the Gaussian kernel with the Hampel loss and $k=5$.

For the KDE and robust KDE, we compute the $L_{\infty}$ influence of $\Y_m$ i.e., $\norminf{f^{n+m}-f^{n}}$ as shown in Figure \ref{fig:influence2} (d). Additionally for each of the $0^{th}$-order and $1^{st}$-order persistence diagrams, we compute the persistence influence of $\Y_m$, i.e., $\Winf\pa{\dgm\pa{f^{n+m}},\dgm\pa{f^{n}}}$ as shown in Figures \ref{fig:influence2} (b, e), and the $1$-Wasserstein influence, i.e., $W_1\pa{\dgm\pa{f^{n+m}},\dgm\pa{f^{n}}}$ as shown in Figures \ref{fig:influence2} (c, f). We refer the reader to Eq. \eqref{eq:p-wasserstein} in Appendix \ref{persistent-homology} for the definition of $W_1$ metric.

For each value of $r$, we generate $100$ such samples and report the average in Figure \ref{fig:influence2}. The results indicate that the robust persistence diagrams, $\dgm\pa{\fn}$, are relatively unperturbed when the outliers are added. It exhibits stability even as $r$ become very large. The KDE persistence diagrams, $\dgm\pa{\barfn}$, on the other hand, are unstable as the outlying noise becomes more extreme.

\begin{figure}[H]
  \centering
  \begin{subfigure}[b]{0.3\linewidth}
    \includegraphics[width=\linewidth]{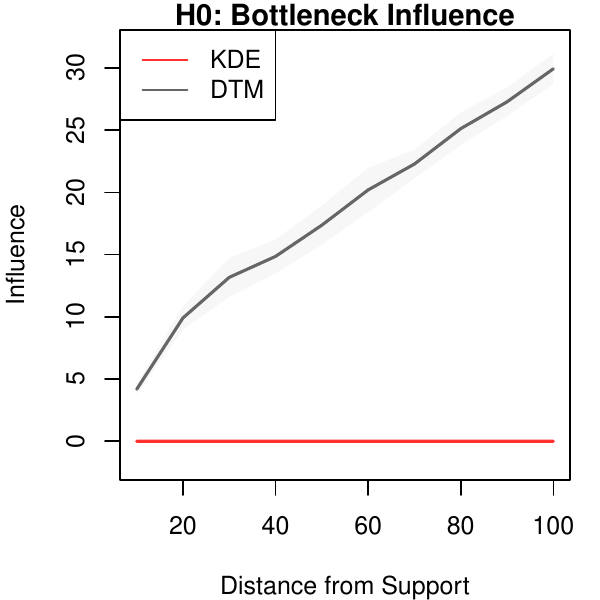}
    \caption{}
  \end{subfigure}
  \begin{subfigure}[b]{0.3\linewidth}
    \includegraphics[width=\linewidth]{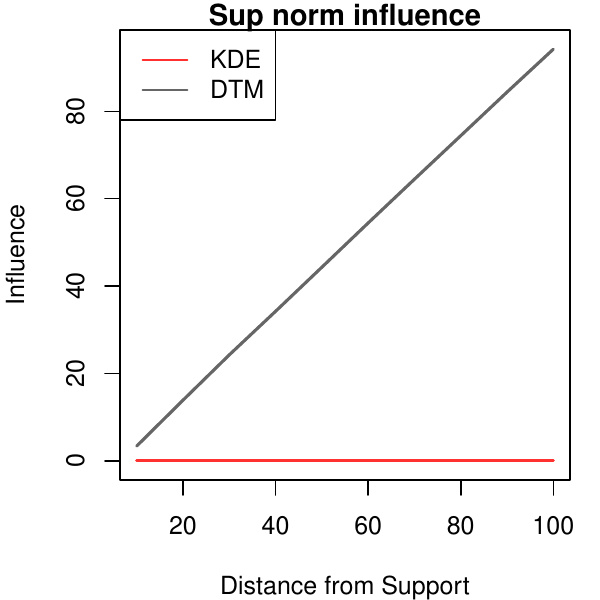}
    \caption{}
  \end{subfigure}
  \begin{subfigure}[b]{0.3\linewidth}
    \includegraphics[width=\linewidth]{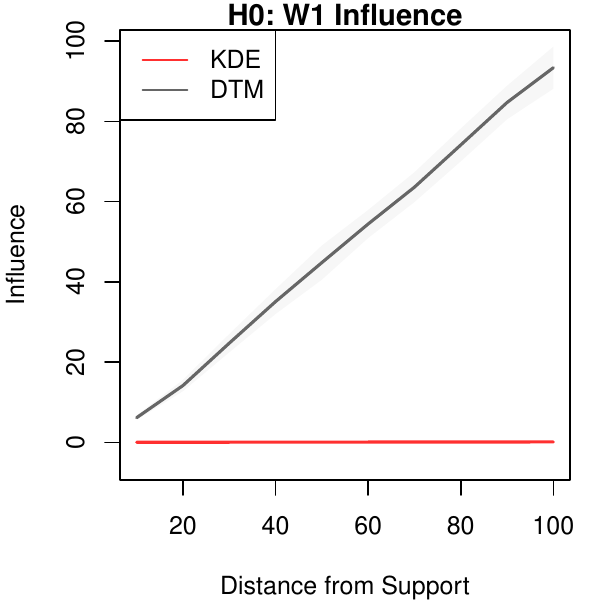}
    \caption{}
  \end{subfigure}
  \begin{subfigure}[b]{0.3\linewidth}
    \includegraphics[width=\linewidth]{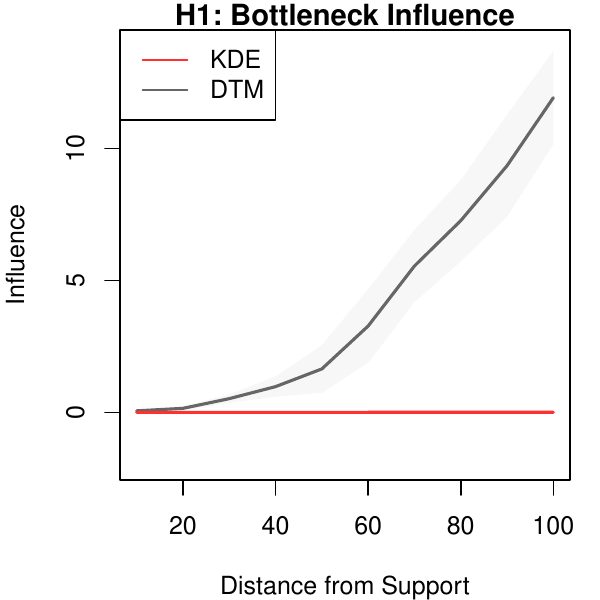}
    \caption{}
  \end{subfigure}
  \begin{subfigure}[b]{0.3\linewidth}
    \includegraphics[width=\linewidth]{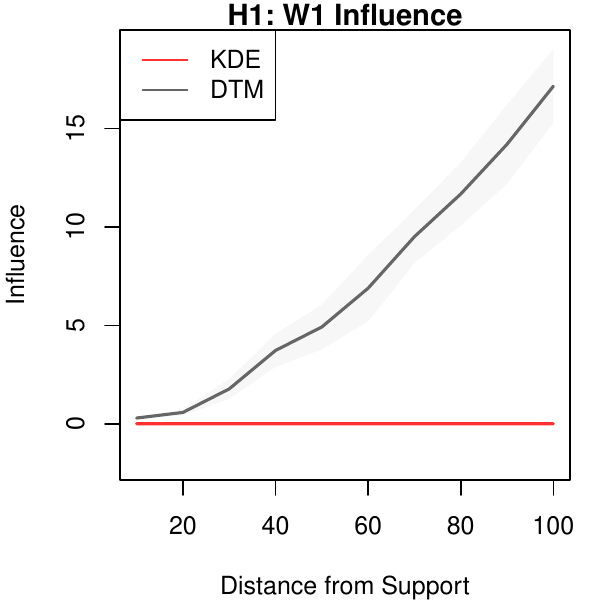}
    \caption{}
  \end{subfigure}
  \caption{For the same data in Figure \ref{fig:influence2}, (a, d) depicts the bottleneck influence for the DTM in contrast to the KDE -- the red line is the same as the one from Figure \ref{fig:influence2} (b, e). Similarly, in (c, e) we see the $W_1$ persistence influence of $\Y_m$ for the DTM in contrast to the KDE. (b) shows the $L_{\infty}$ influence of $\Y_m$ on the DTM. The robust KDE lines were omitted from all plots as it appears to almost merge with the KDE at this scale.}
  \label{fig:dtm-influence}
\end{figure}

As discussed in the Remark~\ref{remark:influence}(iii),
the persistence influence for DTM has a much weaker bound as the outliers become more extreme, and in general is not guaranteed to be bounded. In Figure \ref{fig:dtm-influence} we illustrate the results from the same experiment when the persistence diagrams from DTM is contrasted with the persistence diagrams from the KDE. This analysis is for the same data as that used in Figure \ref{fig:influence}. We remark that even though DTM is highly sensitive to extreme outliers, DTM based filtrations have other remarkable properties, as described in \cite{chazal2017robust}. They are very useful for analyzing persistent homology when one has access to just a single collection of points $\Xn$. For DTM the smoothing parameter is chosen as $m(k) = k/n$ with $k=5$.



\section{Additional Experiments with Robust Persistence Diagrams}
\label{additional-experiments}

In this section, we provide information on some additional experiments with the proposed robust persistence diagrams. The experimental setup is the same as in  Section~\ref{experiments}.

\textbf{Random Circles.} The objective of this simulation is to evaluate the performance of persistence diagrams in a supervised learning task. We select circles $\mathbb{S}_1, \mathbb{S}_2,\dots,\mathbb{S}_{\mathbf{N}}$ randomly in $\R^2$ with centers inside $\pc{0,2}^2$, with the number of such circles, $\mathbf{N}$ uniformly sampled from $\pb{1,2,\dots,5}$. Conditional on $\mathbf{N}=N$, $\Xn$ is sampled uniformly from $\mathbb{S}_1,\dots,\mathbb{S}_N$ with $50\%$ noise in the enclosing square. Two such point clouds are shown in Figure \ref{fig:circles} (a, b). Persistence diagrams $\dgm\pa{\barfn}$ and $\dgm\pa{\fn}$ are constructed for bandwidth $\s(k)$ selected from $k=5,7$, and vectorized in the form of persistence images $\Img\pa{\barfn,h}$, and $\Img\pa{\fn,h}$ for varying bandwidths $h$ \cite{adams2017persistence}. With $\mathbf{N}$ as the response and the persistence images as the input, results from a support vector regression, averaged over 50 random splits, is shown in Figure \ref{fig:circles} (c, d). For a fixed $h$ the robust persistence diagram seems to always contain more predictive information, as observed in the envelope it forms in Figure \ref{fig:circles} (c, d).

{\setlength{\intextsep}{0pt} 
\setlength{\textfloatsep}{0pt} 
\setlength{\abovecaptionskip}{0pt}
\setlength{\belowcaptionskip}{0pt}%
\begin{figure}[H]
  \centering
  \begin{subfigure}[b]{0.24\linewidth}
    \includegraphics[width=\linewidth]{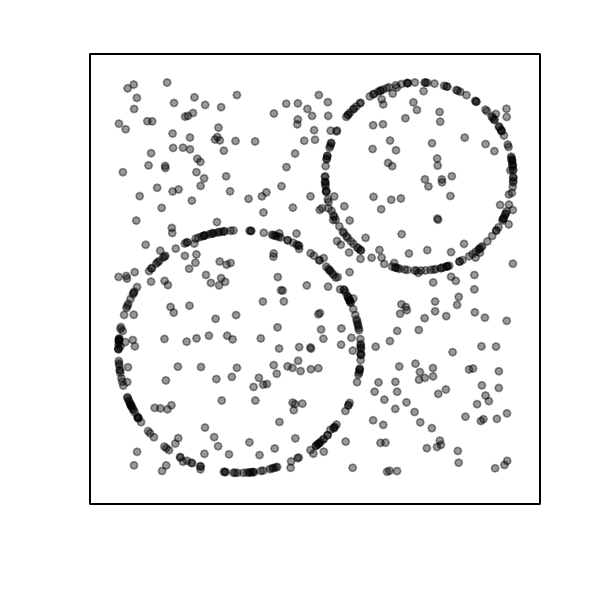}
    \caption{$\Xn$ when $\mathbf{N}=2$}
  \end{subfigure}
  \begin{subfigure}[b]{0.24\linewidth}
    \includegraphics[width=\linewidth]{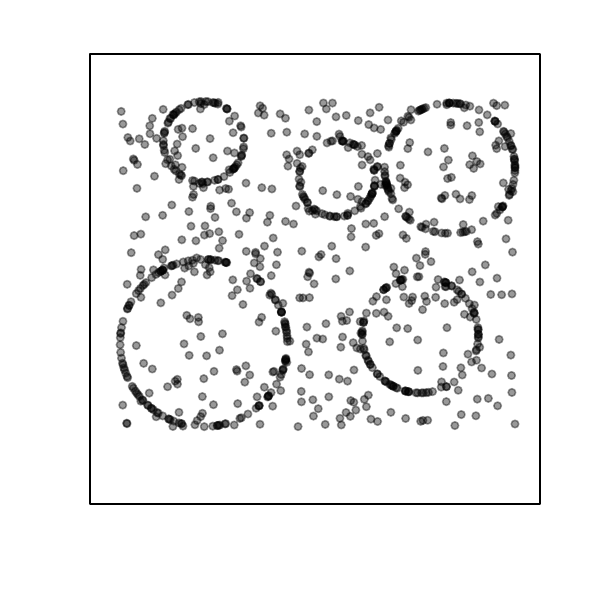}
    \caption{$\Xn$ when $\mathbf{N}=5$}
  \end{subfigure}
  \begin{subfigure}[b]{0.23\linewidth}
    \includegraphics[width=\linewidth]{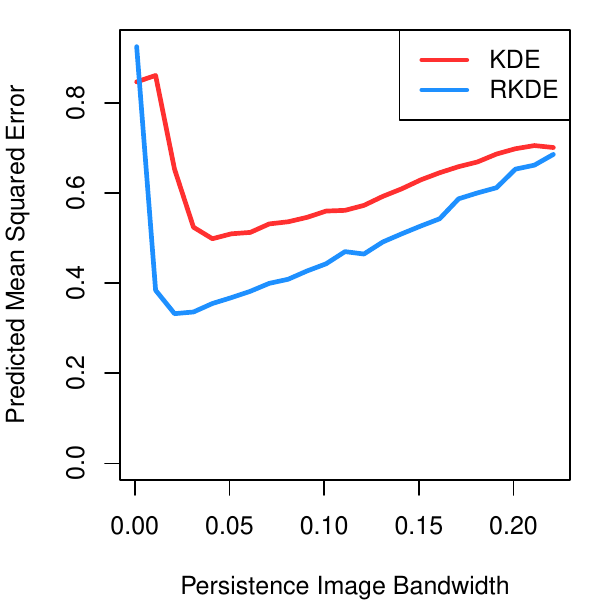}
    \caption{$k=5$}
  \end{subfigure}
  \begin{subfigure}[b]{0.23\linewidth}
    \includegraphics[width=\linewidth]{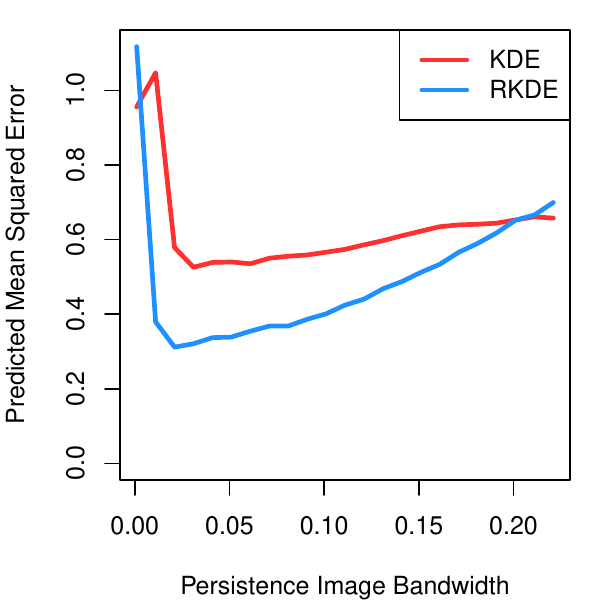}
    \caption{$k=7$}
  \end{subfigure}
\caption{(a, b) A realization $\Xn$ when $\mathbf{N}=2$ and $\mathbf{N}=5$. (c, d) The predicted mean-squared error vs.~the persistence image bandwidth for persistence diagrams in support vector regression.}  \label{fig:circles}
\end{figure}
}


\section{Background on Persistent Homology}
\label{persistent-homology}
Given a set of a points $\Xn = \pb{\xv_1 \dots \xv_n}$ in a metric space $\pa{\X,d}$ their topology is encoded in a geometric object called a simplicial complex $\mathcal{K} \subseteq 2^{\Xn}$.

\begin{samepage}
\begin{definition}
  \textup{\cite{hatcher2005algebraic}}. \ {A \textit{simplicial complex} $\mathcal{K}$ is a collection of simplices $\inner{\sigma}$ i.e. points, lines, triangles, tetrahedra and its higher dimensional analogues, such that
  \begin{enumerate}
    \item $\forall \tau \preccurlyeq \sigma$, $\sigma \in \mathcal{K}$ we have $\tau \in \mathcal{K}$;
    \item $\forall \sigma, \tau \in \mathcal{K}$, we have that $\sigma \cap \tau \preccurlyeq \sigma, \tau$ or $\sigma \cap \tau = \phi$.
  \end{enumerate}
  }
\end{definition}
\end{samepage}

For a given spatial resolution $r>0$, the simplicial complex for $\Xn$, given by $\mathcal{K}\pa{\mathbb{X}_n,r}$, can be constructed in multiple ways. For example, the Vietoris-Rips complex is the simplicial complex
\eq{\mathcal{K}_r = \{\sigma \subseteq \Xn : \bigcap_{\xv \in \sigma} B(\xv,r) \neq \varnothing \},\nn
}
and the \cech{} complex is given by
\eq{
\mathcal{K}_r = \{\sigma \subseteq \Xn : \max\limits_{\xv_i,\xv_j \in \sigma}d\pa{\xv_i,\xv_j} \le r\}.\nn
}
More generally, if $\mathcal{K}$ is a simplicial complex constructed using an approximation of the space $\X$ (e.g., triangulation, surface mesh, grid, etc.), and $\phi:\X\rightarrow \R$ a filter function, $\phi$ induces the map $\phi: \mathcal{K} \rightarrow \R$. Then, $\mathcal{K}_r = \phi\inv\pa{[0,r]}$ encodes the information in the sublevel set of $\phi$ at resolution $r$. Similarly, $\mathcal{K}^r$ encodes the information in the superlevel sets at resolution $r$.

For $0 \le k \le d$, the $k^{th}$-\textit{homology} [\citealp{hatcher2005algebraic}] of a simplicial complex $\mathcal{K}$, given by $H_k\pa{\mathcal{K}}$ is an algebraic object encoding its topology as a vector-space (over a fixed field). Using the Nerve lemma, $H_k\pa{\mathcal{K}\pa{\mathbb{X}_n,r}}$ is isomorphic to the homology of its union of $r$-balls, $H_k\pa{\bigcup_{i=1}^{n}B_{r}\pa{\xv_i}}$. The ordered sequence $\pb{\mathcal{K}\pa{\mathbb{X}_n,r}}_{r > 0}$ forms a \textit{filtration}, encoding the evolution of topological features over a spectrum of resolutions. For $0<r<s$, the simplicial complex $\mathcal{K}\pa{\mathbb{X}_n,r}$ is a \textit{sub-simplicial complex} of $\mathcal{K}\pa{\mathbb{X}_n,s}$. Their homology groups are associated with the inclusion maps
\eq{
\iota_r^s:H_k\pa{\mathcal{K}\pa{\mathbb{X}_n,r}} \hookrightarrow H_k\pa{\mathcal{K}\pa{\mathbb{X}_n,s}},\nn
}
which in turn carry information on the number of non-trivial $k$-cycles. As the resolution $r$ varies, the evolution of the topology is captured in the filtration. Roughly speaking, new cycles (e.g., connected components, loops, voids and higher order analogues) can appear or existing cycles can merge. Formally, a new $k$-cycle  $\sigma_k$ with homology class $\pc{\alpha_k}$ is \textit{born} at $b \in \R$ if $\pc{\alpha_k} \notin \text{Im}(\iota^k_{b-\epsilon,b})$ for all $\epsilon > 0$ and $\pc{\alpha_k} \in \text{Im}(\iota^k_{b,b+\delta})$ for some $\delta > 0$. The same $k$-cycle born at $b$ dies at $d > b$ if $\iota^k_{b,d-\delta}\pa{\pc{\alpha_k}} \notin \text{Im}(\iota^k_{b-\epsilon,d-\delta})$ and $\iota^k_{b,d}\pa{\pc{\alpha_k}} \in \text{Im}(\iota^k_{b-\epsilon,d})$ for all $\epsilon > 0$ and $0 < \delta < d-b$. Persistent homology, $PH_*(\phi)$, is an algebraic module which tracks the persistence pairs $(b,d)$ of births $b$ and deaths $d$ across the entire filtration. By collecting all persistence pairs $(b,d)$, the persistent homology is represented as a persistence diagram
\eq{
\dgm\pa{\mathcal{K}\pa{\Xn}} \defeq \pb{(b,d) \in \R^2 : 0 \le b < d \le \infty}.\nn
}
The persistence diagram is a multiset of points on the space ${\Omega = \pb{(x,y): 0 \le x < y \le \infty}}$, such that each point $(x,y)$ in the persistence diagram corresponds to a distinct topological feature which existed in $\mathcal{K}(\Xn,r)$ for $x \le r < y$. Given a persistence diagram $\mathbf{D}$ and $1 \le p \le \infty$ the \textit{degree-$p$ total persistence} of $\mathbf{D}$ is given by
\eq{
\textup{pers}_p(\mathbf{D}) = \pa{\sum_{(b,d) \in \mathbf{D}}\abs{d-b}^p}^{\f{1}{p}}.\nn
}
The space of persistence diagrams, given by $\mathcal{D}_p = \pb{\mathbf{D} : \textup{pers}_p(\mathbf{D}) < \infty}$, is endowed with the family of $p$-Wasserstein metrics $W_p$. Given two persistence diagrams $\mathbf{D}_1,\mathbf{D}_2 \in \mathcal{D}_p$, the $p$-Wasserstein distance is given by
\eq{
W_p\pa{\mathbf{D}_1,\mathbf{D}_2} \defeq \pa{\inf_{\gamma \in \Gamma}\sum_{\zv \in \mathbf{D}_1 \cup \Delta}\norminf{\zv-\gamma(\zv)}^p}^{\f{1}{p}},
\label{eq:p-wasserstein}
}
where $\Gamma = \pb{\gamma : \mathbf{D}_1 \cup \Delta \rightarrow \mathbf{D}_2 \cup \Delta}$ is the set of all bijections from $\mathbf{D}_1$ to $\mathbf{D}_2$ including the diagonal $\Delta = \pb{(x,y) \in \R^2 : 0 \le x=y \le \infty}$ with infinite multiplicity.

\subsection{Weighted Rips Filtrations}

\begin{figure}
  \centering
  \begin{subfigure}[b]{0.31\linewidth}
    \includegraphics[width=\linewidth]{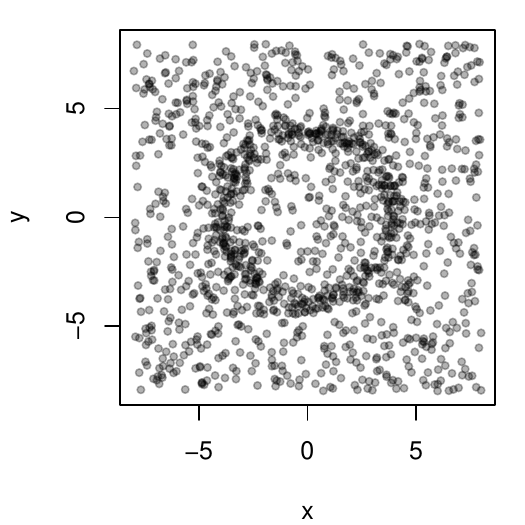}
    \caption{Sample points $\Xn$}
  \end{subfigure}
  \begin{subfigure}[b]{0.31\linewidth}
    \includegraphics[width=\linewidth]{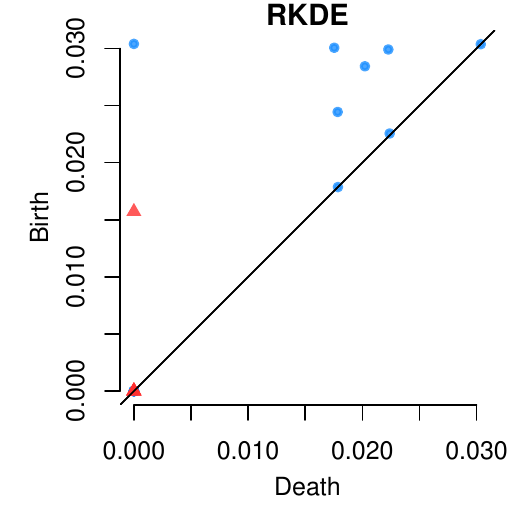}
    \caption{$\dgm\pa{\fn}$}
  \end{subfigure}
  \begin{subfigure}[b]{0.31\linewidth}
    \includegraphics[width=\linewidth]{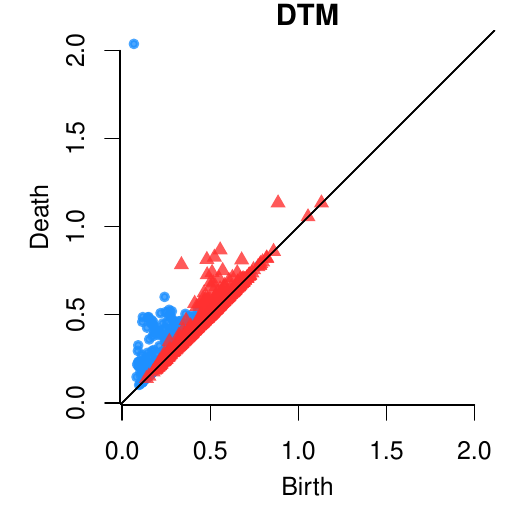}
    \caption{$\dgm\pa{d_{m,n}}$}
  \end{subfigure}
  \par\bigskip
  \begin{subfigure}[b]{0.31\linewidth}
    \includegraphics[width=\linewidth]{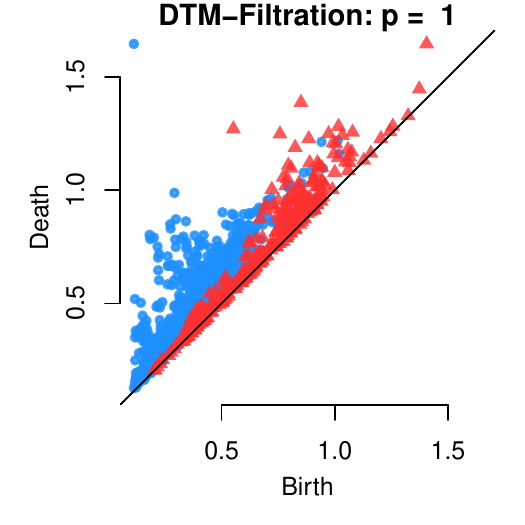}
    \caption{$\dgm\pa{\Xn; d_{n,m}, p=1}$}
  \end{subfigure}%
  \begin{subfigure}[b]{0.31\linewidth}
    \includegraphics[width=\linewidth]{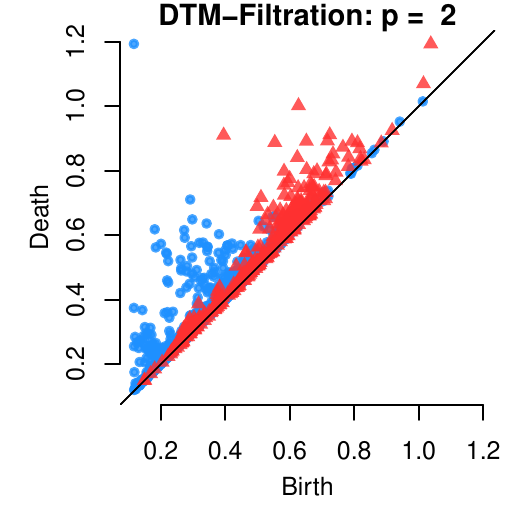}
    \caption{$\dgm\pa{\Xn; d_{m,n}, p=2}$}
  \end{subfigure}
  \begin{subfigure}[b]{0.31\linewidth}
    \includegraphics[width=\linewidth]{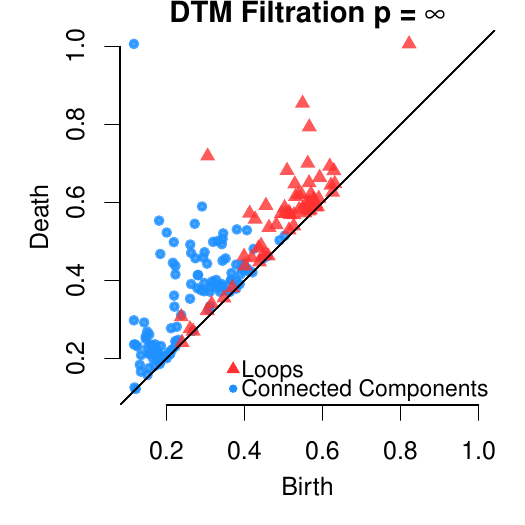}
    \caption{$\dgm\pa{\Xn; d_{m,n}, p=\infty}$}
  \end{subfigure}
  \caption{Points are sampled from a circular region with adverse outlying noise in the enclosing region. The persistence diagrams from sublevel and superlevel filtration from $\fn$ and $d_{n,m}$ respectively are compared with those from the DTM-filtration for $p \in \pb{1,2,\infty}$. The connected components are shown in \textcolor{dblue}{$\bullet$} and loops in \textcolor{red}{$\blacktriangle$}.}
  \label{fig:dtm-filt}
\end{figure}

\begin{figure}
  \centering
  \begin{subfigure}[b]{0.32\linewidth}
    \includegraphics[width=\linewidth]{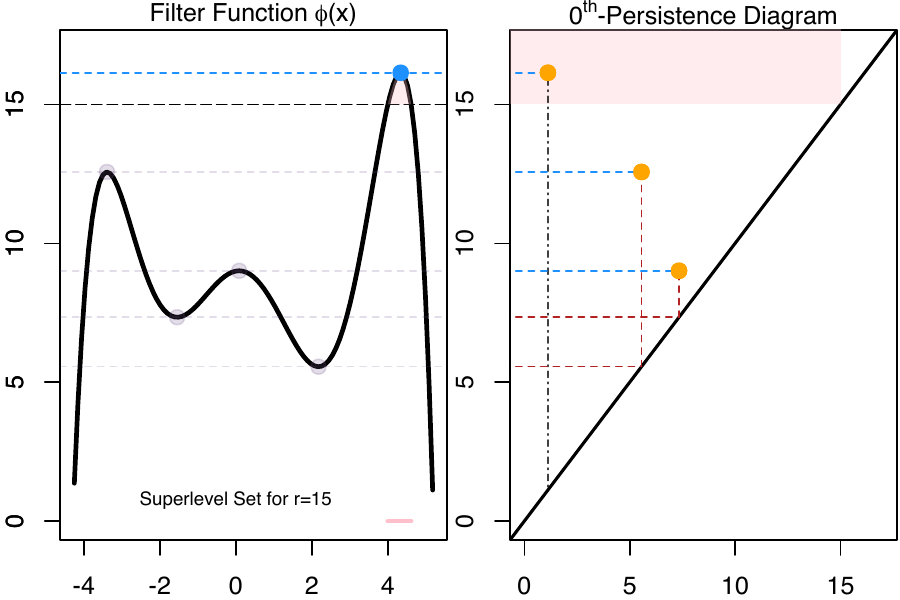}
    \caption{Birth at level $r \approx 15$}
  \end{subfigure}
  \begin{subfigure}[b]{0.32\linewidth}
    \includegraphics[width=\linewidth]{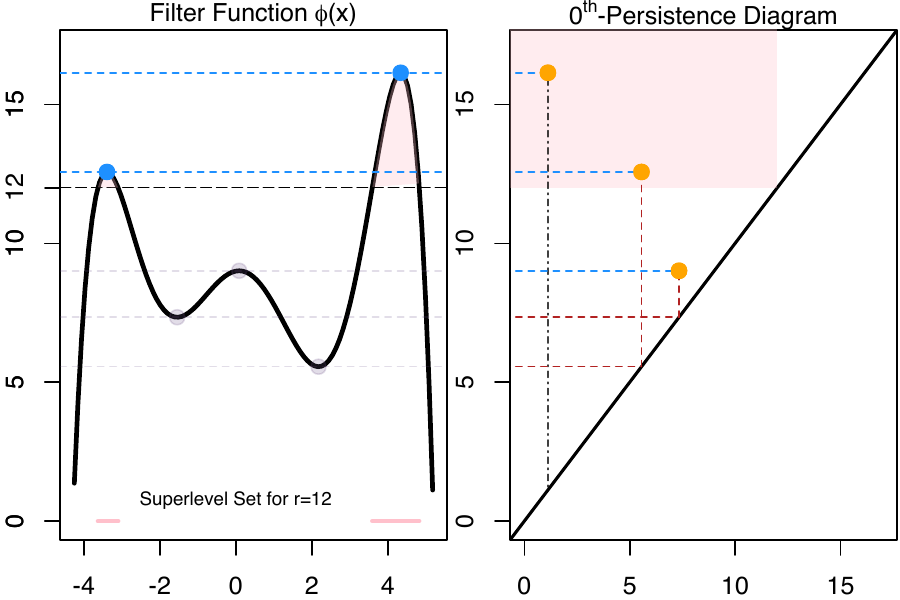}
    \caption{Birth at level  $r \approx 12$}
  \end{subfigure}
  \begin{subfigure}[b]{0.32\linewidth}
    \includegraphics[width=\linewidth]{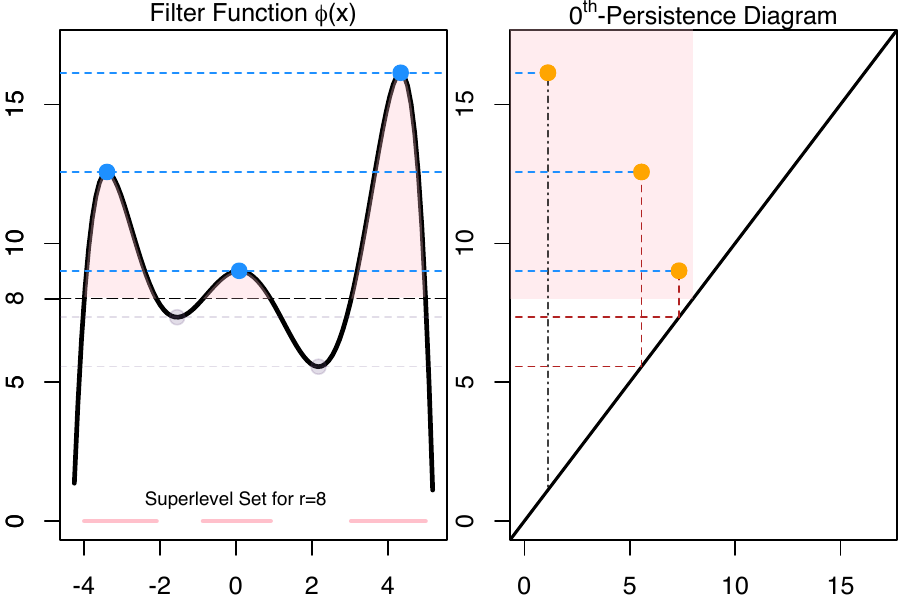}
    \caption{Birth at level  $r \approx 8$}
  \end{subfigure}
  \par\bigskip
  \begin{subfigure}[b]{0.32\linewidth}
    \includegraphics[width=\linewidth]{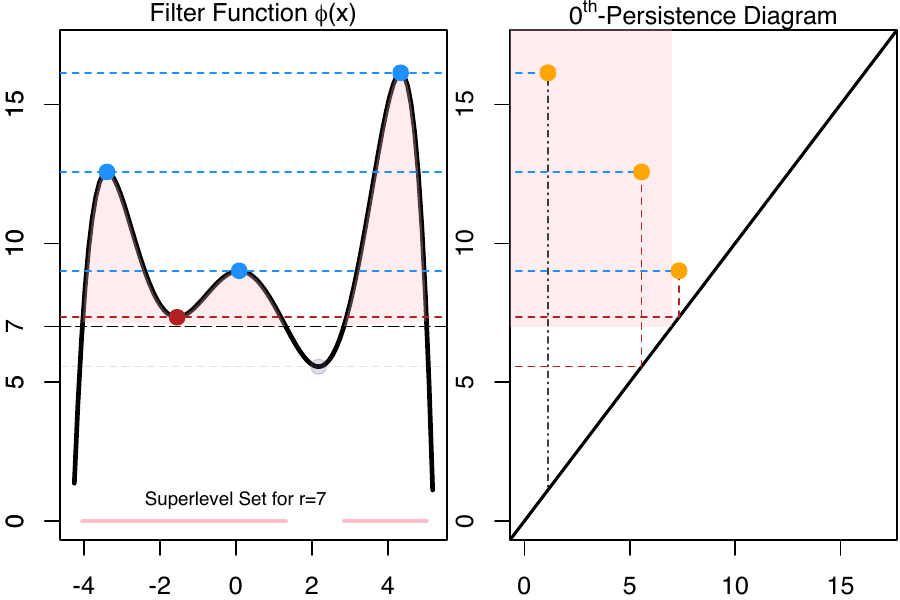}
    \caption{Death at level $r \approx 7$}
  \end{subfigure}%
  \begin{subfigure}[b]{0.32\linewidth}
    \includegraphics[width=\linewidth]{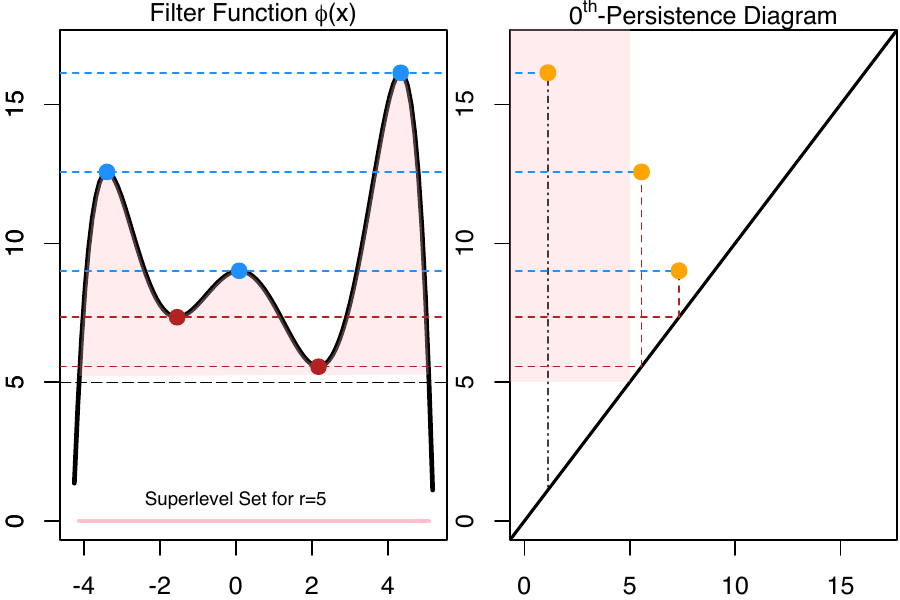}
    \caption{Death at level $r \approx 5$}
  \end{subfigure}
  \begin{subfigure}[b]{0.32\linewidth}
    \includegraphics[width=\linewidth]{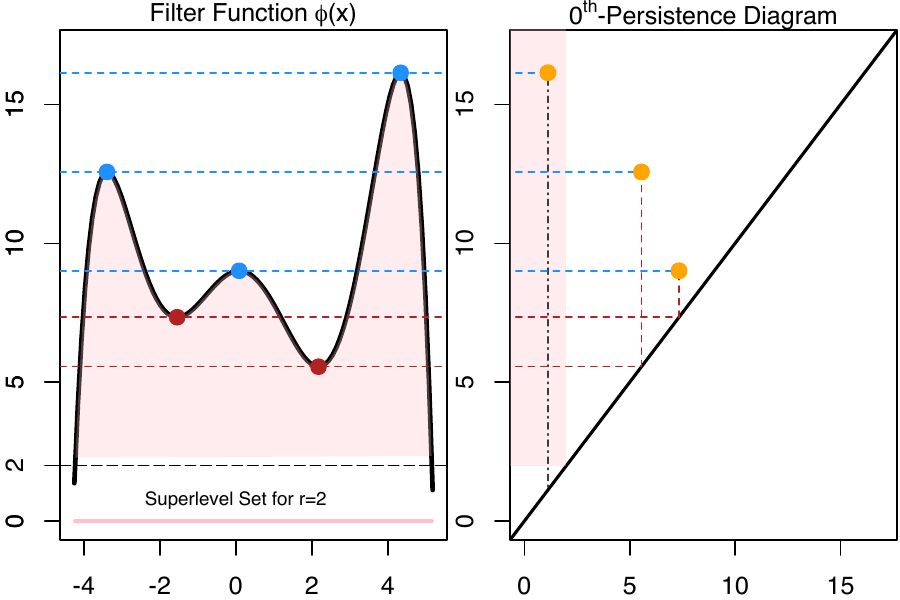}
    \caption{Connected component continues to $r = 0$}
  \end{subfigure}
  \caption{An example for the superlevel filtration of $\phi : \R \rightarrow \R$. (a) As the superlevel set enters $r\approx 15$, the first connected component is born, corresponding to the blue dot on the highest peak of $\phi$. The superlevel set for $r=15$ is depicted in pink below. This is recorded as a birth in the corresponding orange dot enclosed in the pink shaded region of the persistence diagram. (b) As the $r$ enters $r\approx 12$, another connected component is born. This is recorded as the second orange dot in the shaded region of the persistence diagram. (c) Again, at $r \approx 8$, a third connected component is born at the lowest peak of $\phi$. The three connected components in the superlevel set are shaded in pink below the function. The persistence diagram has three orange dots corresponding to these three connected components. (d) As $r$ enters the first valley of $\phi$, depicted by the red dot, two connected components merge (i.e., one of the existing connected components die). By convention, the most recent persistent feature is merged into the older one, i.e., the connected component from (c) merges into the one from (b), and thus, it dies at this resolution. In the persistence diagram, this is noted by the fact that the orange dot born in (c) dies at resolution $r \approx 7$. At this stage, there are only two orange dots in the pink shaded region of the persistence diagram, corresponding to the two pink connected components in the superlevel set of $\phi$. (e) When $r$ enters the second valley of $\phi$, the connected component from (b) merges into the connected component from (a), and form a single connected component. The orange dot in the persistence diagram records the death of this feature. (f) The single connected component persists from then on, and eventually dies at $r=0$. }
  \label{fig:filtration-illustration}
\end{figure}

{For $p \ge 1$ and a weight function $w:\R^d \rightarrow \R$,  the $p^{\text{th}}$-power distance from $\xv \in \Xn$ at resolution $t > 0$ is given by $r_{\xv,w,p}(t) \defeq \pa{t^p - w(\xv)^p}^{\f{1}{p}}$. \citet{anai2019dtm} introduce the weighted-Rips filtration, where the weighted-Rips complex at resolution $t > 0$ is the simplicial complex
\eq{
\mathcal{K}_{t,w,p} \defeq \pb{\s \subseteq \Xn : \bigcap_{\xv \in \sigma} B\pa{\xv, r_{\xv,w,p}(t)} \neq \varnothing}.
}
The weighted-Rips filtration, $\pb{\mathcal{K}_{t,w,p}}_{0 \le t < \infty}$ is used to construct the persistence diagram $\dgm\pa{\Xn;w,p}$. On the computational front, the construction of $\dgm\pa{\Xn;w,p}$ does not depend on the dimension of the underlying space. As a result, weighted-Rips filtrations are very appealing for applications in high dimensions. In addition, the weighted-Rips filtrations obtained by using the distance-to-measure (DTM) as the weight function, i.e., $\dgm\pa{\Xn;d_{m,n},p}$, have some very appealing approximation properties \cite[Theorems 15 \& 20]{anai2019dtm}.

We highlight some key differences between our approach and that in \citet{anai2019dtm}. First, as remarked in \cite[Section 5]{anai2019dtm}, many of the favourable properties of the DTM-filtrations follow from the stability of DTM w.r.t. the Wasserstein Distance. However, it should be noted that stability is inherently different from robustness, as we have described in our analysis using the \textit{persistence influence} in Section~\ref{robustness}, and particularly, in Proposition~\ref{prop:dtm-influence}. In this context, Figure~\ref{fig:dtm-filt} demonstrates the advantage of our proposed approach in the presence of adverse noise. Second, we note that that our implementations of the robust persistence diagrams use the superlevel filtrations of the robust KDE $\fn$ (e.g., see Figure~\ref{fig:filtration-illustration}), in contrast to weighted-Rips filtrations. While latter is arguably better in higher dimensions, it becomes infeasible for large sample sizes. Notwithstanding, the contributions of \cite{anai2019dtm} provide an interesting direction to pursue using the tools presented here to develop \textit{efficient and robust} persistence diagrams.}

\end{document}